\documentclass[12pt]{article}

\usepackage[english]{babel}
\usepackage[T1]{fontenc}
\usepackage{amsmath}
\usepackage{amsfonts}
\usepackage{amssymb}
\usepackage[utf8]{inputenc}
\usepackage{graphicx,color,url}
\usepackage[margin=2.4cm]{geometry}
\usepackage{booktabs}
\usepackage{hyperref}
\usepackage{subfigure}
\usepackage{caption}
\usepackage{tikz}
\usepackage{natbib}

\def\no{\noindent}

\newtheorem{theo}{Theorem}
\newtheorem{rem}{Remark}

\def\calH{\mathcal{H}}

\def\mca{\mathcal{A}}
\def\mcq{\mathcal{Q}}
\def\mco{\mathcal{O}}

\def\pmatrix{\left(\begin{array}}
\def\endpmatrix{\end{array}\right)}
\def\dd{\mathrm{d}}

\def\sech{\mathrm{sech}}

\def\mom{\mu_m}
\def\mon{\mu_n}

\begin{document}
\title{Simple bespoke preservation of two conservation laws}
\author{G.\,Frasca-Caccia\,\quad P.\,E.\,Hydon\,\\[.5cm]
\small
School of Mathematics, Statistics and Actuarial Science\\
\small
University of Kent, Canterbury, CT2 7FS\\
\small}

\maketitle

\begin{abstract} {Conservation laws are among the most fundamental geometric properties of a partial differential equation (PDE), but few known finite difference methods preserve more than one conservation law. All conservation laws belong to the kernel of the Euler operator, an observation that was first used recently to construct approximations symbolically that preserve two conservation laws of a given PDE. However, the complexity of the symbolic computations has limited the effectiveness of this approach. The current paper introduces some key simplifications that make the symbolic-numeric approach feasible. To illustrate the simplified approach, we derive bespoke finite difference schemes that preserve two discrete conservation laws for the Korteweg-de Vries (KdV) equation and for a nonlinear heat equation. Numerical tests show that these schemes are robust and highly accurate compared to others in the literature.}
		
\medskip
\no{\bf Keywords:} 	{Finite difference methods; discrete conservation laws; KdV equation; nonlinear heat equation; porous medium equation.}
	%

\end{abstract}

\section{Introduction}\label{intro}

The main goal of geometric integration is to reproduce, in a numerical approximation, key geometric properties of a given continuous differential problem \citep[see][]{HLW06,BuddPiggott}.

For instance, Hamiltonian ordinary differential equations (ODEs) occur in applications from nano-scale molecular dynamics to the macro-scale of celestial mechanics \citep[see][]{HLW06,BI}. They have two fundamental features: symplecticity of the flow in phase space and constancy of the Hamiltonian function on solutions. Consequently, geometric integration of Hamiltonian ODEs has followed two main approaches, preserving symplecticity and energy respectively.
Symplectic methods are obtained by requiring that the discrete map associated with a given numerical method is symplectic \citep[][]{LeimkuhlerReich,HLW06,Feng,SanzSerna,SanzCalvo}.
Energy conservation has been achieved by using discrete line integral methods \citep[see][]{BI,BIT10,BIT12,BIT15}, time finite element methods \citep[][]{BetschSteinmann,TangChen,TangSun} and discrete gradient methods \citep[][]{DahlbyOwren,Gonzales,McLQuisRobi} such as the Average Vector Field method \citep[see][]{CMcLMcLOQW09,QuispelMcLaren,Ha10}.

These structure-preserving approaches have been extended to Hamiltonian partial differential equations (PDEs)
\citep[][]{BridgesReich,Br97,LeimkuhlerReich}. A particularly powerful approach uses a multisymplectic reformulation of the equations \citep[][]{Bridges,IsSc04,Br97,BridgesReich,LeimkuhlerReich,IKS01, CQT02,AschMcLac,AschMcLac2,QinSun,SunQin}. Alternatively, the method of lines is used to create a semidiscretization, and the resulting Hamiltonian ODEs (in time) are integrated by a symplectic method \citep[][]{Bridges,Cano06,OWW04,QinZhang90,AschMcLac,LuSchmid,BFG13,GJLT09} or energy-conserving method \citep[][]{BBFI,BBFI2,BBFI18,BI,BFI15,BFI152,BFI153,DahlbyOwren,GFC,GX15,Furihata,KoideFur}. For Hamiltonian PDEs, \citet{McLQ} made the useful observation that  \textit{``if the semidiscretization has a semidiscrete energy conservation law, then a discrete gradient method applied to this semidiscretization will have a fully discrete energy conservation law''}.

The benefits of preserving global invariants have been examined for several Hamiltonian PDEs in \citet{FrutosSS,DuranLM,DuranSS} and, in a more general context, in \citet{DuranSSgen}.

The current paper introduces a simple bespoke approach to constructing finite difference schemes that preserve multiple conservation laws of a given PDE. Conservation laws are among the most fundamental features of the PDE, as their origin is topological. Our approach is a simplification of the symbolic-numeric strategy introduced in \citet{GrantPHD} and \citet{GrantHydon} and developed in \citet{Grant}.

There are three advantages to this approach. First, it does not require the PDE to have any special structure, so it is suitable for discretizing PDEs independently of whether or not they possess other geometric structures. Second, the discretizations obtained by using this strategy exactly preserve {\em local} discrete conservation laws. Conserving local features of the continuous PDE gives, in general, a stricter constraint than preserving the corresponding global features. Given suitable boundary conditions, the preservation of local conservation laws also ensures the conservation of the corresponding global invariants. Finally, our approach can be used to seek methods that preserve any number of conservation laws. However, imposing the preservation of more than two conservation laws can considerably increase the complexity of the scheme. For this reason, in this paper, we deal only with methods that preserve two conservation laws, as a reasonable compromise between reliability and complexity of the schemes.

In Section~\ref{method} we review Grant's symbolic-numeric approach and introduce the simplifications that we will use to construct new conservative finite difference schemes. A different strategy, the multiplier method, has been proposed in \citet{Bihlo} to construct conservative finite difference methods for ODEs and PDEs. We briefly discuss the two different approaches. In Section~\ref{KdVsec}, the simplified symbolic-numeric approach is applied to the Korteweg-de Vries (KdV) equation. Several new schemes are constructed and numerical tests are presented to show their effectiveness by comparison with some known methods that preserve only one conservation law. In Section~\ref{heat}, we consider a nonlinear heat equation as an example of a non-Hamiltonian PDE having two conservation laws. A family of two-parameter methods preserving both conservation laws is introduced. (These are easily extended to the more general porous medium equation.) At the end of the section, we present numerical tests that show the conservative properties of the new schemes and comparisons with a standard second-order finite difference method. Some concluding remarks are given in Section~\ref{conclusions}.

\section{How to preserve multiple conservation laws}\label{method} We begin this section with some basic results on conservation laws of partial differential equations (PDEs). After reviewing the general symbolic-numeric strategy for preserving multiple conservation laws of a given scalar PDE, we introduce some simplifications that enable accurate schemes to be derived efficiently. 

We restrict attention to scalar PDEs with two independent variables; the approach generalizes to more variables, but more simplifications may be needed to make the symbolic computations tractable. Consider a PDE for $u(x,t)$,
\begin{equation}\label{PDE}
\mca(x,t,[u])=0,
\end{equation}
where $[u]$ denotes $u$ and finitely many of its derivatives. More generally, square brackets around a differentiable expression denote the expression and finitely many of its derivatives. To simplify the exposition, we will assume that $\mca$ is at most quadratic in $[u]$; the generalization to PDEs that are polynomial in $[u]$ is obvious. For an application of our approach to a PDE that is cubic in $[u]$, see \citet{FCmKdV}.

A conservation law of \eqref{PDE} is a divergence expression,
$$\mbox{Div}\,\mathbf{F} \equiv D_x\{F(x,t,[u])\}+D_t\{G(x,t,[u])\},$$
which is zero on all solutions of (\ref{PDE}); that is,
\begin{equation}\label{CL}
\mbox{Div}\,\mathbf{F} =0\quad \mbox{when}\quad [\mca=0].
\end{equation}
Here $D_x$ and $D_t$ are the total derivatives with respect to $x$ and $t$ respectively:
\begin{eqnarray*}
	D_x&\equiv &\frac{\partial}{\partial x}+u_x\frac{\partial}{\partial u}+u_{xx}\frac{\partial}{\partial u_x}+u_{xt}\frac{\partial}{\partial u_t}+\cdots,\\
	D_t&\equiv &\frac{\partial}{\partial t}+u_t\frac{\partial}{\partial u}+u_{tt}\frac{\partial}{\partial u_t}+u_{tx}\frac{\partial}{\partial u_x}+\cdots.
\end{eqnarray*}
The components $F$ and $G$ are commonly referred to as the {\em flux} and {\em density} respectively. A conservation law (\ref{CL}) is {\em trivial of the first kind} if $F$ and $G$ are zero on solutions of (\ref{PDE}). It is {\em trivial of the second kind} if the divergence in (\ref{CL}) is identically zero without any reference to the PDE. A conservation law is {\em trivial} if and only if it is a linear superposition of the two types of trivial conservation laws. Two conservation laws are {\em equivalent} if they differ by a trivial conservation law. If the conservation law \eqref{CL} amounts to 
\begin{equation}\label{QD}
\mbox{Div}\,\mathbf{F}=\mcq\mca,
\end{equation}
it is said to be in {\em characteristic form} and the multiplier $\mcq$ is called a {\em characteristic} of the conservation law. 
\begin{rem}\label{rempdekov}
	If the PDE is in Kovalevskaya form, integrating any of its conservation laws by parts yields an equivalent conservation law in characteristic form \citep[see][]{olver}. A characteristic, $\mcq$, is trivial if it vanishes on solutions of (\ref{PDE}); two characteristics are equivalent if they differ by a trivial characteristic. If (\ref{PDE}) is in Kovalevskaya form, there is a one-to-one correspondence between equivalence classes of characteristics and equivalence classes of conservation laws \citep[see][]{alonso,olver}. Therefore, characteristics can be used to test the equivalence of conservation laws. 
\end{rem}
A crucial result, for our purposes, is the characterization of the kernel of the Euler operator,
\begin{equation*}
\mathcal E=\sum_{i,j}(-D_x)^i(-D_t)^j\frac{\partial}{\partial u_{x^it^j}}\,,\qquad\text{where}\quad u_{x^it^j}=D_x^{\,i}D_t^{\,j}(u),
\end{equation*}
as the space of total divergences. Consequently, if $\mcq$ is a function such that $$\mathcal E(\mcq\mca)\equiv 0,$$ then there exists $\mathbf{F}$ such that $\mcq\mca=\mbox{Div}\,\mathbf{F}$ and, therefore, $\mcq$ is the characteristic of the corresponding conservation law.

Conservation laws, being defined as divergences (\ref{CL}), are local features of (\ref{PDE}). Their integrals over the spatial domain yield quantities that are globally conserved on solutions (provided that (\ref{PDE}) is coupled with suitable boundary conditions). However, although the local preservation of (\ref{CL}) implies the preservation of the globally conserved quantities, the converse is not true. For this reason, we seek finite difference schemes that preserve discrete analogues of continuous local conservation laws.

For simplicity, we consider only uniform discretizations of the PDE (\ref{PDE}). Relative to a generic lattice point $\mathbf{n}=(m,n)$, the grid points are
$$x_i=x(m+i)=x(m)+i\Delta x,\qquad t_j=t(n+j)=t(n)+j\Delta t,$$
and the approximated values of the dependent variable $u\in \mathbb{R}$ at these points are
\begin{equation*}
u_{i,j}\approx u(x_i,t_j),\quad i,j\in\mathbb{Z}.
\end{equation*}
The forward shift operators $S_m$ and $S_n$ are defined on the lattice by
\[
S_m:(m,n)\mapsto (m+1,n),\qquad S_n:(m,n)\mapsto (m,n+1);
\]
their action extends naturally to $x_i,t_j$ and $u_{i,j}$ as follows:
\[
S_m:(x_i,t_j,u_{i,j})\mapsto (x_{i+1},t_j,u_{i+1,j}),\qquad
S_n:(x_i,t_j,u_{i,j})\mapsto (x_i,t_{j+1},u_{i,j+1}).
\]
Combining $S_m$ with the identity operator,
\[
I:(m,n,x_i,t_j,u_{i,j})\to(m,n,x_i,t_j,u_{i,j}),
\]
yields the forward difference, $D_m$, and the forward average, $\mom$, defined for all functions $f$ by
\[
D_m(f)=\tfrac{1}{\Delta x}\,(S_m-I)(f),\quad\ \mom(f)=\tfrac{1}{2}(S_m+I)(f).
\]
To obtain backward versions of the above operators, compose each with $S_m^{-1}$. Similarly,
\[
D_n(f)=\tfrac{1}{\Delta t}\,(S_n-I)(f),\quad \mon(f)=\tfrac{1}{2}(S_n+I)(f).
\]
All of these operators commute with one another.

Discretizing (\ref{PDE}) by means of a suitable finite difference approximation for the derivatives of the dependent variable, one obtains a partial difference equation (P$\Delta$E),
\begin{equation}\label{PdE}
\widetilde\mca(m,n,[u])=0.
\end{equation}
Here $[u]$ denotes $u_{0,0}$ and a finite number of its shifts; more generally, square brackets around a difference expression denote the expression and finitely many of its shifts.

We seek schemes with the following finite difference analogue of each preserved conservation law:
\begin{equation}\label{dCL}
\mbox{Div}\,\widetilde{\mathbf{F}}\equiv D_m\left(\widetilde{F}(m,n,[u])\right)+ D_n\left(\widetilde{G}(m,n,[u])\right)=0\quad\mbox{when}\quad [\,\widetilde\mca=0\,],
\end{equation}
where tildes represent discretizations of the corresponding continuous terms. The functions $\widetilde{F}$ and $\widetilde{G}$ are respectively the flux and the density of the conservation law (\ref{dCL}).

Just as in the continuous case, a conservation law of (\ref{PdE}) is trivial of the first kind if $\widetilde{F}$ and $\widetilde{G}$ vanish on solutions of (\ref{PdE}) and trivial of the second kind if (\ref{dCL}) is identically satisfied without any reference to (\ref{PdE}) and its shifts \citep[see][]{Hydonbook}. A difference conservation law is trivial if and only if it is a linear combination of trivial conservation laws of these two kinds. Two conservation laws are equivalent if they differ by a trivial conservation law.

A conservation law of (\ref{PdE}) is in characteristic form if
$$\mbox{Div}\,\widetilde{\mathbf{F}}=\widetilde\mcq(m,n,[u])\widetilde\mca.$$
Here $\widetilde\mcq$ is the characteristic, which is trivial if it is zero on all solutions of (\ref{PdE}); two characteristics are equivalent if their difference is a trivial characteristic. 

\begin{rem}\label{GH}\citep{GrantHydon, Hydonbook} P$\Delta$Es that can be solved for a highest shift in one direction (such as explicit P$\Delta$Es) admit a one-to-one correspondence between equivalence classes of characteristics and equivalence classes of conservation laws. Therefore, characteristics can be used to test equivalence for conservation laws of such P$\Delta$Es.
\end{rem}

The key result that underpins the symbolic-numeric approach is due to \citet{Kupershmidt}: similarly to the continuous case, the set of all divergence expressions (\ref{dCL}) over $\mathbb{Z}^2$ is precisely the kernel of the difference Euler operator,
\begin{equation}\label{disEuler}
\mathsf{E}\equiv\sum_{i,j}S_m^{-i}S_n^{-j}\frac{\partial}{\partial u_{i,j}}\,.
\end{equation}
\citep[See][for the generalisation of this result.]{HydonMans} Thus, if a function $\widetilde\mcq$ satisfies $$\mathsf{E}(\widetilde\mcq\widetilde\mca)\equiv 0,$$ there exists $\widetilde{\mathbf{F}}$ such that $\widetilde\mcq \widetilde\mca=\mbox{Div}\,\widetilde{\mathbf{F}}$; therefore $\widetilde\mcq$ is the characteristic of this difference conservation law. For consistency, restrict attention to discretizations $\widetilde{\mcq}$ of $\mcq$; then the difference conservation law $\widetilde{\mcq}\widetilde{\mca}$ is automatically a discretization of the continuous conservation law $\mcq\mca$.

Grant's basic symbolic-numeric approach is straightforward. Choose a stencil of points and consider the most general discretizations on the stencil, $\widetilde\mca $ of the PDE and $\widetilde\mcq$ of the characteristic of the desired conservation law. If the stencil is large enough, there will be some free parameters in the discretizations. To preserve the conservation law, impose the condition $\mathsf{E}(\widetilde{\mcq}\widetilde{\mca})=0$. This condition amounts to a system of algebraic equations that express constraints on the parameters. 
The procedure can be iterated for multiple characteristics, $\mcq_l$, provided that the corresponding system of algebraic equations admits a solution.
Finally, consistency conditions are applied to ensure that  $\widetilde{\mca}$ converges to $\mca$ and each $\widetilde{\mcq}_l$ converges to $\mcq_l$ as the stepsizes $\Delta t$ and $\Delta x$ tend to zero; these give further constraints on the free parameters. In this way, bespoke finite difference schemes for a given PDE may be derived by symbolic computation.

In more detail, the basic method is as follows. Having chosen a stencil, the most general discretizations of the PDE (\ref{PDE}) and the characteristics are based on Taylor series expansions of the grid function about the point $(x(m),t(n))\equiv(x_0,t_0)$:
\begin{align}\nonumber
u_{i,j}&\approx u(x_i,t_j)=u(x_0+i\Delta x,t_0+j\Delta t)\\\label{Taylor}
&=u+i\Delta x\,u_x+j\Delta t\,u_t+\frac{(i\Delta x)^2}{2!}u_{xx}+\frac{(j\Delta t)^2}{2!}u_{tt}+i\Delta x\,j\Delta t\,u_{xt}+\cdots\Big\vert_{(x_0,t_0)}.
\end{align}
For a rectangular stencil of points defined by $i=A,\ldots, B$ and $j=C,\ldots, D$, linear terms in $\mca$ and $\mcq$ are approximated by linear combinations, with undetermined coefficients, of terms of the form (\ref{Taylor}):
\begin{equation}\label{linapp}
\frac{\partial^{r+s}}{\partial x^r\partial t^s}u\approx\frac{1}{\Delta x^r}\frac{1}{\Delta t^s}\sum_{i=A}^B\sum_{j=C}^D\alpha_{i,j}u_{i,j},
\end{equation}
where the coefficients $\alpha_{i,j}$ depend on $r$ and $s$.

If quadratic terms appear in $\mca$ or in $\mcq$ (as happens in our examples), we need to look at products of Taylor expansions:
\begin{align}\label{Taylor2}\nonumber
u_{i,j}u_{k,l}=&\,u^2+(i+k)\Delta x\,uu_x+(j+l)\Delta t\,uu_t+ik\Delta x^2u_x^2+(jk+il)\Delta x\Delta t\,u_xu_t+jl\Delta t^2u_t^2\\
&+\frac{\Delta x^2(i^2+k^2)}{2!}uu_{xx}+(ij+kl)\Delta x\,\Delta t\,uu_{xt}+\frac{\Delta t^2(j^2+l^2)}{2!}uu_{tt}+\cdots\Big\vert_{(x_0,t_0)}.
\end{align}
Just as for the linear terms, quadratic quantities in $\mca$ and $\mcq$ are replaced by linear combinations, with undetermined coefficients, of terms of the form (\ref{Taylor2}):
\begin{equation}\label{quadapp}
\frac{\partial^{m+n}u}{\partial x^m\partial t^n}\frac{\partial^{r+s}u}{\partial x^r\partial t^s}\approx\frac{1}{\Delta x^{m+r}\Delta t^{n+s}}\sum_{j=C}^D\sum_{i=A}^B\left(\sum_{k=i}^B\beta_{i,j,k,j}u_{i,j}u_{k,j}+\sum_{l=j+1}^D\sum_{k=A}^B\beta_{i,j,k,l}u_{i,j}u_{k,l}\right),
\end{equation}
with the coefficients $\beta_{i,j,k,l}$ depending on $r$ and $s$.

For the right hand sides of (\ref{linapp}) and (\ref{quadapp}) to approximate the corresponding left hand sides, we also need to impose a number of consistency conditions on the coefficients $\alpha_{i,j}$ and $\beta_{i,j,k,l}$. 

\begin{rem}
	Terms involving higher powers of $u_{i,j}$ may be added, provided that they vanish as $\Delta x$ and $\Delta t$ tend to zero. For simplicity, such terms are not included here.
\end{rem}

Having set $\widetilde{\mcq}$ and $\widetilde{\mca}$ to be the discretizations of $\mcq$ and $\mca$ respectively, one must solve
\begin{equation}\label{E0}
\mathsf{E}(\widetilde{\mcq}\widetilde{\mca})=0,
\end{equation}
where $\mathsf{E}$ is the difference Euler operator (\ref{disEuler}). In general, this is not easy. Typically, even for a PDE (\ref{PDE}) that is only quadratic in $[u]$, (\ref{E0}) amounts to a large system of nonlinear algebraic equations \citep[see][]{GrantPHD,Grant}. In principle, such systems can be solved by finding a Groebner basis \citep{BuchbergerKauers,BuchbergerKauers2,CoxLittleOShea,Mansfield}. However, the calculation of the Groebner basis may take a huge amount of memory and a very long computation time. As the use of this approach is limited mainly by the cost of the symbolic computation, it is helpful to impose some additional assumptions. For instance, \citet{Grant} describes some symmetry-based ans\"{a}tze that can simplify the Groebner basis calculation. 

There is one other approach to constructing finite difference schemes that can preserve multiple conservation laws for systems of PDEs. This is the multiplier method, introduced in \citet{Bihlo}, which is as follows. Given a scalar PDE,
\begin{equation}\label{PDEmult}
\mca =0,
\end{equation}
that has one conservation law in the form
\begin{equation}\label{CLmult}
D_x F+ D_t G=\mcq\mca,
\end{equation}
let $\widetilde{\mcq}$, $\widetilde{F}$ and $\widetilde{G}$ be finite difference approximations of $\mcq$, $F$ and $G$, respectively. Then, provided that $\widetilde{\mcq}^{-1}$ exists on the whole domain of definition of (\ref{PDEmult}),
\begin{equation}\label{multmeth}
\widetilde{\mca}\equiv\widetilde{\mcq}^{-1}\left\{D_m\left(\widetilde{F}\right)+D_n\left(\widetilde{G}\right)\right\}=0
\end{equation}
is a finite difference approximation of (\ref{PDEmult}) that preserves the conservation law (\ref{CLmult}).
This method has been applied to find conservative schemes for the inviscid Burgers' equation and momentum-preserving schemes for the KdV equation. It has also been applied to systems of $r$ PDEs, such as the two-dimensional shallow-water equations, to preserve $s$ conservation laws, with $s\leq r$.

The multiplier method has the advantage of being simple to use.  However, there are two main disadvantages. First, it cannot find schemes preserving $s$ conservation laws for systems of $r$ PDEs with $s>r$; in particular, it cannot preserve multiple conservation laws for a scalar PDE. Second, it requires the characteristic to be nonzero throughout the domain. For characteristics that involve the dependent variable, one cannot identify points where the characteristic is zero \textit{a priori}.

By contrast, Grant's approach is in principle able to find all conservative finite difference methods on the chosen stencil, provided that one is able to solve (\ref{E0}). The procedure can be iterated to select discretizations that preserve further conservation laws, provided that the corresponding condition (\ref{E0}) admits a solution for each characteristic.
Moreover, there is no need to choose a particular and arbitrary discretization of densities and fluxes, as these can be reconstructed from the characteristics \citep[see][]{Hydon}.

\begin{figure}
	\center{
		\begin{tikzpicture}
		\draw[fill] (0,0) circle [radius=0.075];
		\draw[fill] (2,0) circle [radius=0.075];
		\draw[fill] (4,0) circle [radius=0.075];
		\draw[fill] (6,0) circle [radius=0.075];
		\draw[fill] (8,0) circle [radius=0.075];
		\draw[fill] (0,3) circle [radius=0.075];
		\draw[fill] (0,4.5) circle [radius=0.075];
		\draw[fill] (2,4.5) circle [radius=0.075];
		\draw[fill] (4,4.5) circle [radius=0.075];
		\draw[fill] (6,4.5) circle [radius=0.075];
		\draw[fill] (8,4.5) circle [radius=0.075];
		\draw[fill] (0,1.5) circle [radius=0.075];
		\draw[fill] (2,1.5) circle [radius=0.075];
		\draw[fill] (4,1.5) circle [radius=0.075];
		\draw[fill] (6,1.5) circle [radius=0.075];
		\draw[fill] (8,1.5) circle [radius=0.075];
		\draw[fill] (2,3) circle [radius=0.075];
		\draw[fill] (4,3) circle [radius=0.075];
		\draw[fill] (6,3) circle [radius=0.075];
		\draw[fill] (8,3) circle [radius=0.075];
		\draw[thick,dotted] (0,0)--(8,0);
		\draw[thick,dotted] (0,1.5)--(8,1.5);
		\draw[thick,dotted] (0,3)--(8,3);
		\draw[thick,dotted] (0,4.5)--(8,4.5);
		\draw[thick,dotted] (0,0)--(0,4.5);
		\draw[thick,dotted] (2,0)--(2,4.5);
		\draw[thick,dotted] (4,0)--(4,4.5);
		\draw[thick,dotted] (6,0)--(6,4.5);
		\draw[thick,dotted] (8,0)--(8,4.5);
		\node [below right] at (3.9,1.45) {\scriptsize{$(x,t\!-\!{\Delta t/2})$}};
		\node [below right] at (3.9,2.2) {\scriptsize{$(x,t)$}};
		\node [below] at (3,2.2) {\scriptsize{$(x\!-\!{\Delta x/2},t)$}};
		\draw (4,2.25) circle [radius=0.125];
		\draw (3.85,1.35)--(4.15,1.65);
		\draw (3.85,1.65)--(4.15,1.35);
		\draw (2.85,2.10)--(3.15,2.4);
		\draw (2.85,2.4)--(3.15,2.1);
		\end{tikzpicture}
		\caption{Example of a rectangular stencil. Conservation laws are preserved to second-order accuracy at the central point $(x, t)$ (circle), densities and fluxes respectively at $(x, t-\Delta t/2) $ and $(x-\Delta x/2, t)$ (crosses).}
		\label{stencil}
	}
\end{figure}
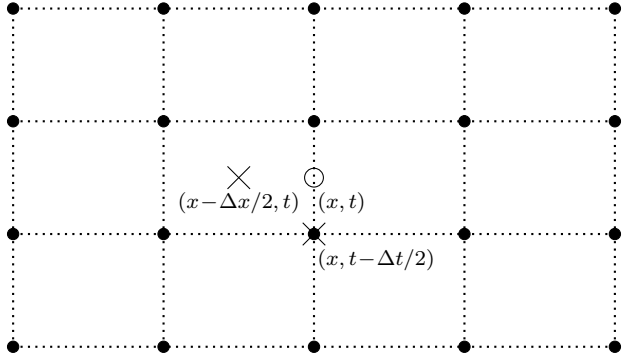

To simplify Grant's approach, we adopt a strategy that reduces the number of variables and the computational cost of solving the system of nonlinear equations. This is achieved by first looking for second-order accurate approximations only, building in consistency from the outset. If the stencil is as compact as possible, this immediately determines the discretizations of the highest order derivatives. The problem can be further simplified by restricting the approximations of some terms in $\widetilde{\mca}$ and $\widetilde{\mcq}$ to use only points in a sub-stencil that is as compact as possible. In particular, by approximating nonlinear terms using as few points as possible, the number of variables may be considerably reduced to the point of being able to solve (\ref{E0}) with a fast symbolic computation that does not need a Groebner basis.
\begin{rem}\label{remsurf} Conservation laws of a given PDE are evaluated on hypersurfaces. In the discrete case, the smallest ``surface'' on which a conservation law can be evaluated locally is the convex hull of the stencil \citep{Hydonbook}. Therefore, for all conservative schemes that are presented in this paper, the consistency conditions are imposed so as to approximate the conservation laws {(and their characteristics)} to second-order accuracy at the centre of the {rectangular} stencil. 
\end{rem}
\begin{rem}\label{rem2ord}
	Obtaining second-order accurate approximations of the conservation laws at the centre $(x,t)$ of the stencil is equivalent to finding second-order accurate approximations of the corresponding densities and fluxes at the points $(x,t-\Delta t/2)$ and $(x- \Delta x/2,t)$ respectively (see (\ref{dCL})). Figure~\ref{stencil} shows an example of a rectangular stencil. The circle denotes the centre, i.e. the point where we require second order approximations of the PDE and of the characteristics (and hence of the conservation laws). The crosses denote the points where we require second order approximations of the corresponding densities and fluxes. These are not necessarily lattice points.
\end{rem}

In the next {two} sections, we use the above simplifications to derive approximations that preserve two conservation laws of some well-known nonlinear wave equations, and show that these conservative methods can generate robust, highly-accurate schemes.
\section{KdV equation}\label{KdVsec}

In this section we exploit the strategy introduced in Section~\ref{method} to develop conservative schemes for the KdV equation,
\begin{equation}\label{KdV}
\mca\equiv u_t+uu_x+u_{xxx}=0, \quad (x,t)\in\Omega\equiv[a,b]\times[0,\infty).
\end{equation}
These schemes are tested for two benchmark problems; they compare favourably with two well-known schemes that each preserve only one conservation law.

Equation (\ref{KdV}) has an infinite number of conservation laws. The first three, in increasing order, are
\begin{align}\label{KdvCL1}
& D_t(G_1)+D_x(F_1) \equiv D_t(u)+D_x\left(\frac{1}2 u^2 + u_{xx}\right)=0,\\\label{KdVCL2}
& D_t(G_2)+D_x(F_2) \equiv D_t\left(\frac{1}2 u^2\right)+D_x\left(\frac{1}3u^3+uu_{xx}-\frac{1}2u_x^2\right)=0,\\\label{KdVCL3}
& D_t(G_3)+D_x(F_3) \equiv D_t\left(\frac{u^3}3+uu_{xx}\right)+D_x\left(\frac{u^4}4+u_xu_t-uu_{xt}+u^2u_{xx}+u_{xx}^2\right)=0,
\end{align}
which can be written in characteristic form (\ref{QD}) with characteristics
\begin{equation}\label{KdVchar}
\mcq_1=1,\qquad \mcq_2=u,\qquad \mcq_3=u^2+2u_{xx},
\end{equation}
respectively. For a water wave problem, conservation laws (\ref{KdvCL1})-(\ref{KdVCL3}) describe the local conservation of mass, momentum and energy, respectively \citep[see][]{DrazinJohnson}. As these conservation laws have a physical meaning, it seems particularly desirable to preserve them.

When (\ref{KdV}) is coupled with suitable (e.g. periodic or zero) boundary conditions, integrating (\ref{KdvCL1})--(\ref{KdVCL3}) over the spatial domain gives the global conservation of, respectively,
\begin{equation}\label{glinvkdv}
\int G_1\,\dd x=\int u\,\dd x,\qquad \int G_2\,\dd x=\int \frac{1}{2}u^2\,\dd x, \qquad \int G_3\,\dd x=\int \frac{1}{3}u^3+uu_{xx}\, \dd x.
\end{equation}
It is well-known that (\ref{KdV}) possesses the Hamiltonian structure
\begin{equation*}
u_t=D_x\frac{\delta}{\delta u}\calH_1,
\end{equation*}
where ${\delta}/{\delta u}$
is the variational derivative and
\begin{equation}\label{Hamiltonian}
\calH_1=\int-\frac{1}2\left(\frac{1}3u^3-u_x^2\right)\dd x,
\end{equation}
is the Hamiltonian functional. Equivalently, one can use the alternative Hamiltonian functional
\begin{equation}\label{Hamiltonian2}
\calH_1=\int-\frac{1}2\left(\frac{1}3u^3+uu_{xx}\right)\dd x=-\frac{1}{2}\int G_3\,\dd x.
\end{equation}
With this choice of functional, the conservation law (\ref{KdVCL3}) implies the preservation of $\calH_1$.

The KdV equation (\ref{KdV}) can also be written in another Hamiltonian form \citep[see][]{olver},
$$u_t=\left(D_x^3+\frac{2}{3}uD_x+\frac{1}{3}u_x\right)\frac{\delta}{\delta u}\calH_2,$$
with the Hamiltonian
\begin{equation}\label{lowH}
\calH_2=\int-\frac{1}{2}u^2\,\mathrm{d}x=-\int G_2\,\mathrm{d}x.
\end{equation}
The conservation law (\ref{KdVCL2}) implies that $\calH_2$ is preserved.

These are special instances of the following well-known general result. Given a scalar Hamiltonian evolution equation for $u(x,t)$,
\begin{equation}\label{Hamev}
u_t=\mathcal{D}\frac{\delta}{\delta u}\calH,\qquad \calH=\int H([u]) \,\mathrm{d}x,
\end{equation}
where $\mathcal{D}$ is a skew-adjoint differential operator (with respect to the $L^2$ inner product), the Hamiltonian $\calH$ is constant (provided that some technical conditions are satisfied).

To discretize the KdV equation (\ref{KdV}), one first needs to set the stencil. Having done this, we will use our simplified version of Grant's approach to construct two types of scheme: our \textit{energy-conserving schemes} preserve discrete versions of the conservation laws (\ref{KdvCL1}) and (\ref{KdVCL3}), while our \textit{momentum-conserving schemes} preserve discrete versions of (\ref{KdvCL1}) and (\ref{KdVCL2}).

\subsection{Conservative methods for the KdV equation}\label{KdVmeth}

\noindent\textbf{8-point schemes}\\

\begin{figure}
	\center{
		\begin{tikzpicture}
		\draw[fill] (0,0) circle [radius=0.075];
		\draw[fill] (2,0) circle [radius=0.075];
		\draw[fill] (4,0) circle [radius=0.075];
		\draw[fill] (6,0) circle [radius=0.075];
		\draw[fill] (0,2) circle [radius=0.075];
		\draw[fill] (2,2) circle [radius=0.075];
		\draw[fill] (4,2) circle [radius=0.075];
		\draw[fill] (6,2) circle [radius=0.075];
		\node [below] at (3,0) {\tiny{$(-1/2,0)$}};
		\node [below right] at (3,1) {\tiny{$(-1/2,1/2)$}};
		\node [below] at (2,1) {\tiny{$(-1,1/2)$}};
		\draw (3,1) circle [radius=0.175];
		\draw (1.87,0.87)--(2.13,1.13);
		\draw (1.87,1.13)--(2.13,0.87);
		\draw (2.87,-0.13)--(3.13,0.13);
		\draw (2.87,0.13)--(3.13,-0.13);
		\draw (3.87,-0.13) rectangle (4.13,0.13);
		\node [below right] at (4,0) {\tiny{$(0,0)$}};
		
		\end{tikzpicture}
		\caption{The most compact rectangular stencil for (\ref{KdV}). Conservation laws are preserved to second-order accuracy at the central point $(-1/2, 1/2)$ (circle), densities and fluxes respectively at $(-1/2,0) $ and $(-1, 1/2)$ (crosses).}
		\label{stencilcomp}
	}
\end{figure}
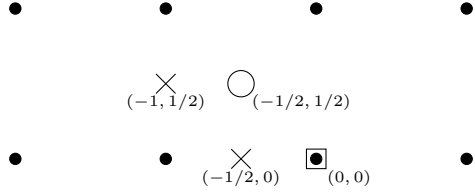

The most compact rectangular stencil consists of 8 points, as shown in Fig.~\ref{stencilcomp}. Here and henceforth, grid points are labelled with respect to the lattice point denoted with a square. From Remark~\ref{rem2ord}, we seek second-order approximations of characteristics, densities and fluxes at $(-1/2,1/2)$, $(-1/2,0)$ and $(-1,1/2)$, respectively.

\noindent\textit{Energy-conserving schemes}\\

To simplify the symbolic computations for energy-conserving schemes on the 8-point stencil, we use the approximations
\[
{u_{xx}}\approx\mu_n\left(D_m^2u_{-2,0}\right),\qquad {u_{xx}}\approx\mu_m\left(\mu_n\left(D_m^2u_{-2,0}\right)\right),
\]
in $\widetilde{F_1}$ and $\widetilde{Q_3}$, respectively. The remaining terms in the approximations of $G_1$ at $(-1/2,0)$, $F_1$ at $(-1,1/2)$
and $\mcq_3$ at $(-1/2,1/2)$ are obtained from (\ref{linapp}) and (\ref{quadapp}) by requiring that the coefficients $\alpha_{i,j}$ and $\beta_{i,j,k,l}$ satisfy all consistency conditions for second-order accuracy. This yields families of approximations that depend on just a few undetermined coefficients.

The discretizations of the conservation laws (\ref{KdvCL1}) and (\ref{KdVCL3}) at $(-1/2,1/2)$ are taken to be of the following form (setting $\widetilde{\mcq_1}=1$):
\[
\widetilde{\mca}=D_m\left(\widetilde{F_1}\right)+D_n\left(\widetilde{G_1}\right)=0,\qquad\qquad
\widetilde{\mcq_3}\widetilde{\mca}=0.
\]
As $\widetilde{\mca}$ is defined to be a discrete conservation law, the condition $\mathsf{E}(\widetilde{\mca})\equiv 0$ holds
for any choice of the remaining coefficients in $\widetilde{G_1}, \widetilde{F_1}$ and $\widetilde{Q_3}$.
We then find these undetermined coefficients by solving
\begin{equation*}
\mathsf{E}(\widetilde{\mcq}_3\widetilde{\mca})\equiv 0.
\end{equation*}
This constraint determines all remaining coefficients. So only one scheme of this form preserves (\ref{KdvCL1}) and (\ref{KdVCL3}) to second-order accuracy at the centre of the stencil:
\begin{equation}\label{E-C8}
\mbox{EC}_8\equiv\widetilde{\mca}=D_m\left(\widetilde{F_1}\right)+D_n\left(\widetilde{G_1}\right)=0,
\end{equation} 
where
\begin{equation}\label{13Grant}
\widetilde{F_1}=\frac{1}{3}\,\mu_n\left(\left(\mu_m^2u_{-2,0}\right)^2\right)+\frac{1}{6}\left(\mu_m^2u_{-2,0}\right)\left(\mu_m^2u_{-2,1}\right)+\mu_nD_m^2u_{-2,0},\qquad \widetilde{G_1}=\mu_mu_{-1,0}.
\end{equation}
The scheme $\mbox{EC}_8$ preserves the following discrete version of the conservation law (\ref{KdVCL3}):
\begin{align}\label{2clKdVmethodgen}
\widetilde{\mcq_3}\widetilde{\mca}=& D_m\widetilde{F_3}+ D_n\widetilde{G_3}=0,
\end{align}
where
\begin{align}\nonumber
\widetilde{\mcq_3}=&\,2\mu_m\widetilde{F_1},\\\nonumber
\widetilde{F_3}=&\,\widetilde{F_1}^2+\left(\mu_nD_m\mu_mu_{-2,0}\right)\left(D_n\mu_m^2u_{-2,0}\right)-\left(\mu_n\mu_m^2u_{-2,0}\right)\left(D_nD_mu_{-2,0}\right)\\\nonumber
&\,+\frac{\Delta x^2}{6}(\mu_n\mu_m^2u_{-2,0})\left\{(2\mu_mD_nu_{-1,0})(\mu_nD_m\mu_mu_{-2,0})-D_n(\mu_mu_{-1,0}D_m\mu_mu_{-2,0})\right\},\\\nonumber
\widetilde{G_3}=&\,\frac{1}{3}{(\mu_mu_{-1,0})\mu_m\left(\left(\mu_m^2u_{-2,0}\right)^2\right)}+(\mu_mu_{-1,0})(\mu_mD_m^2u_{-2,0}).
\end{align}
The last term in the flux $\widetilde{F_3}$ vanishes as the spatial stepsize tends to zero, and does not correspond to an expression in the continuous flux. The scheme $\mbox{EC}_8$ is equivalent to one introduced in \citet{Grant}.

When (\ref{KdV}) is coupled with zero or periodic boundary conditions, the scheme $\mbox{EC}_8$ preserves at each time step the following discretization of the Hamiltonian (\ref{Hamiltonian2}):
\begin{equation}\label{discHGrant}
\widetilde{\calH_1}(j)=-\frac{\Delta x}2\sum_i\left(\frac{1}{3}{(\mu_mu_{i-1,j})\mu_m\left(\left(\mu_m^2u_{i-2,j}\right)^2\right)}+(\mu_mu_{i-1,j})(\mu_mD_m^2u_{i-2,j})\right).
\end{equation}
In general, $\mbox{EC}_8$ fails to preserve the conservation law (\ref{KdVCL2}), even to first order: given any approximation,
$$\widetilde{\mcq_2}=\sum_{i=-2}^1\sum_{j=0}^1\eta_{i,j}u_{i,j},\qquad \widetilde{\mcq_2}= \mcq_2 +\mco(\Delta x,\Delta t),$$
the condition $E(\widetilde{\mcq_2}\widetilde{\mca})=0$ cannot be satisfied when $\widetilde{\mca}$ is given by (\ref{E-C8}).
\bigskip

\noindent\textit{Momentum-conserving schemes}\\

\citet{Grant} introduced several momentum-conserving schemes, including a one-parameter family obtained by using the most compact second-order approximation of $u_t$ in (\ref{KdV}). This amounts to
\begin{equation*}
\mbox{MC}_8(\alpha)\equiv\widetilde{\mca}=D_m\left(\widetilde{F_1}\right)+D_n\left(\widetilde{G_1}\right)=0,
\end{equation*}
where
\begin{align*}\nonumber
\widetilde{F_1}=&\,\frac{1}{6}(\mu_nu_{-1,0})\mu_n(u_{-2,0}+u_{-1,0}+u_{0,0})+\mu_nD_m^2u_{-2,0}\\\nonumber
&\,+\alpha\Delta x^2\left\lbrace(\mu_nu_{0,0})(D_m^2\mu_nu_{-2,0})+D_m((\mu_nu_{-2,0})(D_m\mu_nu_{-2,0}))\right\rbrace,\\
\widetilde{G_1}=&\,\mu_mu_{-1,0}.
\end{align*}
For any value of $\alpha$, these methods preserve the discrete momentum conservation law
$$\widetilde{\mcq_2}\widetilde{\mca}= D_m\left(\widetilde{F_2}\right)+ D_n\left(\widetilde{G_2}\right)=0,$$
with
\begin{align*}
\widetilde{\mcq_2}=&\,\mu_m\mu_nu_{-1,0}\,,\\
\widetilde{F_2}=&\, \frac{1}{3}(\mu_m\mu_nu_{-2,0})(\mu_m\mu_nu_{-1,0})\mu_n(u_{-1,0}+2\alpha\Delta x^2D_m^2u_{-2,0})\\
&\,+(\mu_n\mu_m^2u_{-2,0})(\mu_nD_m^2u_{-2,0})-\frac{1}{2}\mu_m\left((\mu_nD_mu_{-2,0})^2\right),\\
\widetilde{G_2}=&\,\frac{1}{2}(\mu_mu_{-1,0})^2.
\end{align*}
For zero or periodic boundary conditions, these schemes preserve the following discretization of the Hamiltonian (\ref{lowH}) at each time step:
\begin{equation}\label{discHMC8}
\widetilde{\calH_2}(j)=-\frac{\Delta x}2\sum_i\left(\mu_mu_{i-1,j}\right)^2.
\end{equation}
The local truncation error of the scheme $\mbox{MC}_8(\alpha)$ is $\mco(\Delta x^2)+\mco(\Delta t^2)$. Restricting attention to the case $\Delta t \ll\Delta x$, the truncation error can be reduced considerably by choosing $\alpha$ optimally. No choice of $\alpha$ eliminates the second-order terms identically, so the optimal value will depend on the particular problem.\\

\noindent\textbf{10-point schemes}\\

\begin{figure}
\center{
\begin{tikzpicture}
\draw[fill] (0,0) circle [radius=0.075];
\draw[fill] (2,0) circle [radius=0.075];
\draw[fill] (4,0) circle [radius=0.075];
\draw[fill] (6,0) circle [radius=0.075];
\draw[fill] (8,0) circle [radius=0.075];
\draw[fill] (0,2) circle [radius=0.075];
\draw[fill] (2,2) circle [radius=0.075];
\draw[fill] (4,2) circle [radius=0.075];
\draw[fill] (6,2) circle [radius=0.075];
\draw[fill] (8,2) circle [radius=0.075];
\node [below right] at (4,1) {\tiny{$(0,1/2)$}};
\node [below] at (3,1) {\tiny{$(-1/2,1/2)$}};
\draw (4,1) circle [radius=0.175];
\draw (2.87,0.87)--(3.13,1.13);
\draw (2.87,1.13)--(3.13,0.87);
\draw (3.87,-0.13)--(4.13,0.13);
\draw (3.87,0.13)--(4.13,-0.13);
\draw (3.87,-0.13) rectangle (4.13,0.13);
\node [below right] at (4,0) {\tiny{$(0,0)$}};

\end{tikzpicture}
\caption{The 10-point rectangular stencil for (\ref{KdV}). Conservation laws are preserved with higher order at the central point $(0, 1/2)$ (circle), densities and fluxes respectively at $(0,0) $ and $(-1/2, 1/2)$ (crosses).}
\label{stencil10}
}
\end{figure}
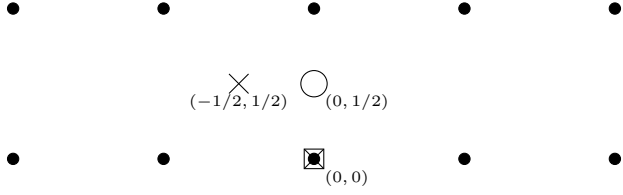

To find new schemes that preserve two conservation laws, one must use a wider stencil. This is beyond what can be tackled in full generality, but the symbolic computations are made tractable (indeed, fast) by the simplifications that we have introduced. Adding one further pair of nodes in the spatial direction gives the 10-point stencil in Fig.~\ref{stencil10}. Hence, according to Remark~\ref{rem2ord}, our goal is to find second-order approximations of characteristics, densities and fluxes at $(0,1/2)$, $(0,0)$ and $(-1/2,1/2)$, respectively.
\bigskip

\noindent\textit{Energy-conserving schemes}\\

Just as for the 8-point schemes, let
\begin{equation*}
\widetilde\mca=D_m\left(\widetilde{F_1}\right)+D_n\left(\widetilde{G_1}\right),
\end{equation*}
so that (\ref{KdvCL1}) is preserved for any choice of the undetermined coefficients. The preservation of (\ref{KdVCL3}), obtained by requiring that 
\begin{equation}\label{eul3}
\mathsf{E}\left(\widetilde{\mcq_3}\widetilde{\mca}\right)\equiv 0,
\end{equation}
is simplified by setting the sub-stencils for $\widetilde{G_1}$ and the quadratic term in $\widetilde{\mcq_3}$ to be as compact as possible, given that the approximations must be second-order:
\begin{align}\label{G110p}
\widetilde{G_1}&=u_{0,0},\\\nonumber
\widetilde{\mcq_3}&=\xi(u_{0,0}^2+u_{0,1}^2)+(1-2\xi)u_{0,0}u_{0,1}+\frac{1}{\Delta x^2}\sum_{i=-2}^2\sum_{j=0}^1\gamma_{i,j}u_{i,j},\quad \xi\in\mathbb{R}.
\end{align}
The undetermined coefficients in $\widetilde{F_1}$ and $\widetilde{\mcq_3}$ are obtained by solving (\ref{eul3}). 
This yields a one-parameter family of schemes,
\begin{equation*}
\widetilde{\mca}(\lambda)=D_m\left(\widetilde{F_1}\right)+D_n\left(\widetilde{G_1}\right)=0,
\end{equation*}
with
\begin{equation*}
\widetilde{F_1}=\mu_m\varphi_{-1,0},\qquad \widetilde{G_1}=u_{0,0},
\end{equation*}
where
$$\varphi_{-1,0}=\frac{u_{-1,1}^2+u_{-1,0}^2+u_{-1,0}u_{-1,1}}{6}+D_m^2\mu_nu_{-2,0}+\lambda D_nD_m\mu_mu_{-2,0},$$
and $\lambda=\mco(\Delta x^2,\Delta t^2)$. These schemes preserve
\begin{align*}
\widetilde{\mcq_3}\widetilde{\mca}=& D_m\widetilde{F_3}+ D_n\widetilde{G_3}=0,
\end{align*}
where
\begin{align*}
\widetilde{\mcq_3}=&\,2\varphi_{0,0}\,,\\\nonumber
\widetilde{F_3}=&\,\varphi_{0,0}\varphi_{-1,0}+(D_m\mu_nu_{-1,0})(D_n\mu_mu_{-1,0})-(\mu_m\mu_nu_{-1,0})(D_nD_mu_{-1,0})\\
&\,+\lambda(D_nu_{0,0})(D_nu_{-1,0}),\\\nonumber
\widetilde{G_3}=&\,\frac{1}{3}u_{0,0}^3+u_{0,0}D_m^2u_{-1,0}\,.
\end{align*}
For zero or periodic boundary conditions, these schemes preserve at each time step
\begin{equation}\label{discH}
\widetilde{\calH_1}(j)=-\frac{\Delta x}2\sum_i\left(\frac{1}3u_{i,j}^3+u_{i,j}\,D_m^2u_{i-1,j}\right),
\end{equation}
but none of them preserves the conservation law (\ref{KdVCL2}).

Assuming for simplicity that $\Delta t\ll\Delta x$, the leading term in the local truncation error amounts to $\mco(\Delta x^2)$. This suggests that by setting $\lambda=\alpha \Delta x^2$, one may be able to remove at least part of this error by choosing $\alpha\in\mathbb{R}$ optimally. However, Taylor expansion shows that no choice of $\alpha$ will give a higher order method. Indeed, the optimal value depends on the initial conditions. In the results section, we write the one-parameter family of schemes as
$$\mbox{EC}_{10}(\alpha)\equiv\widetilde{\mca}(\alpha\Delta x^2).$$

The scheme $\mbox{EC}_{10}(0)$ was originally found by the Discrete Variational Derivative method \citep[see][]{Furihata}. \citet{DahlbyOwren} proved that the Furihata scheme can also be derived by the Average Vector Field method, which approximates \eqref{Hamev} by $D_nu_{0,0}=\widetilde{\mathcal{D}}(\widetilde{\delta}\widehat{\mathcal{H}})$, where the operator $\widetilde{\mathcal{D}}$ is skew-adjoint with respect to the $\ell^2$ inner product. In this case,
\[
\widetilde{\mathcal{D}}=D_m\mu_mS_m^{-1},\qquad \widetilde{\delta}\widehat{\mathcal{H}}=-\frac{1}6\left\{u_{0,1}^2+u_{0,1}u_{0,0}+u_{0,0}^2\right\}-D_m^2\mu_nu_{-1,0}\,.
\]
None of the other $\mbox{EC}_{10}(\alpha)$ schemes can be derived in this way. 
\\

\noindent\textit{Momentum-conserving schemes}\\

To simplify the derivation of momentum-conserving schemes on the 10-point stencil in Fig.~\ref{stencil10}, use the following approximations in which $\widetilde{G_1}$, $\widetilde{\mcq_2}$, and the quadratic term in $\widetilde{F_1}$ are compact\footnote{A family of second-order schemes depending on 8 free parameters can be found by removing the compactness assumption on the quadratic term in $\widetilde{F_1}$.}:
\begin{align}\label{G110mom}
\widetilde{G_1}&\,=u_{0,0},\\\label{Q210p}
\widetilde{\mcq_2}&\,=\mu_nu_{0,0},\\\nonumber
\widetilde{F_1}&\,=\sum_{j=0}^1\sum_{i=-1}^0\sum_{k=i}^0\beta_{i,j,k,j}u_{i,j}u_{k,j}+\sum_{i=-1}^0\sum_{k=-1}^0\beta_{i,0,k,1}u_{i,0}u_{k,1}+\frac{1}{\Delta x^2}\sum_{i=-2}^1\sum_{j=0}^1\alpha_{i,j}u_{i,j},
\end{align}
Proceeding as before, one obtains a two-parameter family of momentum-conserving methods:
\begin{equation*}
\widetilde\mca(\lambda,\nu)=D_m\left(\widetilde{F_1}\right)+D_n\left(\widetilde{G_1}\right)=0,
\end{equation*}
with $\widetilde{G_1}$ as defined in (\ref{G110mom}) and
\[
\widetilde{F_1}=\,\frac{(\mu_nu_{-1,0})^2+(\mu_nu_{0,0})^2+(\mu_nu_{-1,0})(\mu_nu_{0,0})}6+\mu_nD_m^2\mu_mu_{-2,0}
+D_nD_m(\lambda u_{-1,0}+\nu D_m^2u_{-2,0});
\]
here $\lambda=\mco(\Delta x^2,\Delta t^2)$ and $\nu=\mco(\Delta x^2,\Delta t^2)$.
These schemes preserve
\[
\widetilde{\mcq_2}\widetilde{\mca}= D_m\widetilde{F_2}+ D_n\widetilde{G_2}=0,
\]
with $\widetilde{\mcq_2}$ as given in (\ref{Q210p}) and
\begin{align*}
\widetilde{F_2}=&\,\frac{1}{3}(\mu_nu_{0,0})(\mu_nu_{-1,0})(\mu_n\mu_mu_{-1,0})
+\frac{1}2\left\lbrace(\mu_nu_{0,0})(D_m^2\mu_nu_{-2,0})+(\mu_nu_{-1,0})(D_m^2\mu_nu_{-1,0})\right\rbrace\\
&-\frac{1}2(D_m\mu_nu_{-1,0})^2+\lambda\left\lbrace(\mu_m\mu_nu_{-1,0})(D_nD_mu_{-1,0})-\frac{1}{2}D_n(\mu_mu_{-1,0}D_mu_{-1,0})\right\rbrace\\\nonumber
&+\nu\Big\lbrace\!(\mu_nu_{0,0})(D_nD_m^3u_{-2,0})\!-\!(D_m\mu_nu_{-1,0})(D_nD_m^2u_{-1,0})\!\\
&+\!\frac{1}2D_n[D_mu_{-1,0}D_m^2u_{-1,0}-u_{0,0}D_m^3u_{-2,0}]\Big\rbrace,\\
\widetilde{G_2}=&\,\frac{1}2u_{0,0}^2+\frac{1}{2}u_{0,0}D_m^2(\lambda u_{-1,0}+\nu D_m^2u_{-2,0}).
\end{align*}

For suitable boundary conditions, the $\widetilde{\mca}(\lambda,\nu)$ scheme preserves at each time step
\begin{equation}\label{discHMC10}
\widetilde{\calH_2}(j)=-\frac{\Delta x}2\sum_i\left({u_{i,j}^2}+u_{i,j}D_m^2(\lambda u_{i-1,j}+\nu D_m^2u_{i-2,j})\right).
\end{equation}
To simplify the main sources of local truncation error, we will consider only $\Delta t\ll \Delta x$; then the leading term is $\mco(\Delta x^2)$, so we define the two-parameter family of schemes
$$\mbox{MC}_{10}(\alpha,\beta)\equiv\widetilde\mca(\alpha\Delta x^2,\beta\Delta x^2).$$
Again, it is not possible to obtain higher order methods for any choice of the free parameters; the optimal values depend on the particular problem.

It turns out that the $\mbox{MC}_{10}(0,0)$ scheme can also be derived by the Average Vector Field Method, approximating the right-hand side of \eqref{Hamev} by
\[
\widetilde{\mathcal{D}}\widetilde{\delta}\widehat{\mathcal{H}}=D_m^3S_m^{-2}\mu_m(\widetilde{\delta}\widehat{\mathcal{H}})+\frac{2}3S_m^{-1}\mu_m\left\{(\mu_m\mu_nu_{0,0})D_m(\widetilde{\delta}\widehat{\mathcal{H}})\right\}+\frac{1}3D_m(\mu_m\mu_nu_{-1,0})\,\widetilde{\delta}\widehat{\mathcal{H}},
\]
where
\[\widetilde{\delta}\widehat{\mathcal{H}}=-\mu_nu_{0,0}.
\]
In this case, the skew-adjoint operator $\widetilde{\mathcal{D}}$ is not constant.

\subsection{Numerical tests}

In this subsection, two benchmark solutions are used to illustrate the effectiveness of the schemes developed in Section~\ref{KdVmeth} by comparison with two well-known schemes that each preserve a single conservation law. These are the multisymplectic scheme proposed in \citet{AschMcLac2,AschMcLac}, which we rewrite as 
\begin{equation}\label{Asch}
D_m\left(\frac{1}{2}\mu_m(\mu_m\mu_nu_{-2,0})^2+D_m^2\mu_nu_{-2,0}\right)+D_n\left(\mu_m^3 u_{-2,0}\right)=0,
\end{equation}
and the narrow box scheme, defined in the same references, which amounts to
\begin{equation}\label{narr}
D_m\left(\frac{1}{2}(\mu_n u_{-1,0})^2+D_m^2\mu_nu_{-2,0}\right)+D_n\left(\mu_m u_{-1,0}\right)=0.
\end{equation}
Both the multisymplectic scheme (\ref{Asch}) and the narrow box scheme (\ref{narr}) are defined on the 8-point stencil in Fig.~\ref{stencilcomp}, and preserve a discrete version of the mass conservation law (\ref{KdvCL1}).

Each scheme considered in this section is solved by using the Newton method, simplified by using a ``frozen'' Jacobian. This procedure is computationally attractive because the inversion of the Jacobian is performed just once for a single instance of the iterative method. The iterations are run until the error reaches full machine accuracy (up to rounding errors) in double precision. For each of our numerical experiments, the computational cost is approximately the same for all of the schemes.

In the following, we consider (\ref{KdV}) subject to periodic boundary conditions. We evaluate the error in the solution at the final time $t=T$ as
\begin{equation}\label{errkdvsol}
\left.\frac{\|u-u_{\text{exact}}\|}{\|u_{\text{exact}}\|}\right\vert_{t=T}.
\end{equation}
For a grid with $M$ points in space and $N$ points in time, the errors in the invariants (\ref{glinvkdv}) are
\begin{equation}\label{errkdvclCLP}
\text{Err}_\ell=\Delta x\max_{j=1,\ldots,N}\left|\sum_{i=1}^M \left(\widetilde{G_\ell}(x_i,t_j)-\widetilde{G_\ell}(x_i,t_1)\right)\right|,\qquad \ell=1,2,3.
\end{equation}
Where some of the discrete densities $\widetilde{G_1}$, $\widetilde{G_2}$ and $\widetilde{G_3}$ are undefined because the considered scheme does not preserve the corresponding conservation laws, we instead evaluate the corresponding errors as
\begin{align}\nonumber
\text{Err}_1=&\,\Delta x\max_{j=1,\ldots,N}\left|\sum_{i=1}^M (v_{i,j}-v_{i,1})\right|,\\\label{errkdvNCLP}
\text{Err}_2=&\,\Delta x\max_{j=1,\ldots,N}\left|\frac{1}2\sum_{i=1}^M  (v_{i,j}^2-v_{i,1}^2)\right|,\\\nonumber
\text{Err}_3=&\,\Delta x\max_{j=1,\ldots,N}\left|\sum_{i=1}^M \left( \frac{v_{i,j}^3}3+v_{i,j}D_m^2(v_{i-1,j})-\frac{v_{i,1}^3}3-v_{i,1}D_m^2(v_{i-1,1})\right)\right|.
\end{align}
Here $v_{i,j}=u_{i,j}$ for the 10-point stencil in Fig.~\ref{stencil10} and $v_{i,j}=\mu_mu_{i-1,j}$ for the 8-point stencil in Fig.~\ref{stencilcomp}, where $u_{i,j}\simeq u(a+i\Delta x, j\Delta t)$; subscripts denote shifts from the point $(x,t)=(a,0)$.
Note that $\text{Err}_2$ shows how well each scheme preserves the corresponding discretization of the Hamiltonian $\calH_2$, because
$$\max_{j} |\widetilde{\calH_2}(j)-\widetilde{\calH_2}(1)|=\text{Err}_2.$$
Similarly, $\text{Err}_3$ shows how well $\calH_1$ is preserved, because
$$\max_{j} |\widetilde{\calH_1}(j)-\widetilde{\calH_1}(1)|=\text{Err}_3/2.$$

As a first numerical test, we consider equation (\ref{KdV}) for $t\in[0,2]$, with periodic boundary conditions over the interval $[-20,20]$ and the initial condition
\begin{equation}\label{initsolKdV}
u(x,0)=3c\,\sech^2\left(\frac{\sqrt c}{2}(x+d)\right).
\end{equation}
The exact solution of (\ref{KdV}) with initial condition (\ref{initsolKdV}) on an infinite domain is a single soliton,
\begin{equation}\label{ex1solKdV}
u_{\text{exact}}(x,t)=3c\,\sech^2\left(\frac{\sqrt c}{2}(x-ct+d)\right).
\end{equation}

Each scheme is solved for the parameters $c=d=5$ with stepsizes $\Delta x=0.1$ and $\Delta t=0.01$. 
For this problem, the the solution errors for $\mbox{MC}_8(\alpha_1)$, $\mbox{MC}_{10}(\alpha_2,\beta_2)$ and $\mbox{EC}_{10}(\alpha_3)$ are minimised when $\alpha_1=-0.069$, $(\alpha_2,\beta_2)=(0.39,0.04)$ and $\alpha_3=0.12$. When the solution error cannot be evaluated because the exact solution is unknown, another criterion is needed to optimise the free parameters. For instance, to minimize the error in the non-preserved conservation law, the approximate parameter values are $\alpha_1=-0.073$, $(\alpha_2,\beta_2)=(0.21,0.03)$ and $\alpha_3=0.17$.
	
\begin{table}[t!]
\caption{Errors in conservation laws and solution at the final time $T=2$, when solving {\rm(\ref{KdV})} with initial condition {\rm(\ref{initsolKdV})} (for $c=d=5$) and periodic boundary conditions over $[-20,20]$, using various schemes with $\Delta x=0.1$, $\Delta t=0.01$. The non-zero values of the free parameter minimize the error indicated by a star.}
	\label{comp1solkdv}	
\small
	\centerline{\begin{tabular}{ccccc}
	\hline
			Method &  $\text{Err}_1$ & $\text{Err}_2$ & $\text{Err}_3$  &  Solution error\\ 
\hline
			$\mbox{EC}_8$	& 9.24e-14  & 0.0019 & 6.71e-12 &  0.0964\\
			$\mbox{MC}_8(0)$	&  1.24e-13   & 6.82e-13 & 0.0308 & 0.0584\\
			$\mbox{MC}_8(-0.069)$	& 1.24e-13  & 1.05e-12 & 0.0028 &  \phantom{$^*$}0.0052$^*$ \\
			$\mbox{MC}_8(-0.073)$	& 8.17e-14 & 6.25e-13 & \phantom{$^*$}0.0014$^*$ & 0.0063\\
			Furihata; $\mbox{EC}_{10}(0)$	& 8.17e-14   & 9.93e-04 & 2.16e-12 &  0.0217\\
			$\mbox{EC}_{10}(0.12)$	& 7.46e-14  & 3.40e-04 & 2.61e-12 &  \phantom{$^*$}0.0020$^*$\\
			$\mbox{EC}_{10}(0.17)$	& 7.46e-14 &\phantom{$^*$}7.29e-05$^*$ & 2.27e-12 & 0.0095\\
			$\mbox{MC}_{10}(0,0)$	& 6.75e-14  & 3.69e-13 & 0.0033 &  0.0335\\
			$\mbox{MC}_{10}(0.39,0.04)$	& 7.46e-14  & 3.41e-13 & 0.0127 &  \phantom{$^*$}0.0033$^*$\\
			$\mbox{MC}_{10}(0.21,0.03)$	& 6.39e-14 &4.26e-13 & \phantom{$^*$}6.07e-04$^*$ & 0.0224\\
			Multisymplectic & 1.24e-13  & 7.03e-04 & 0.0436 &  0.0385\\
			Narrow box & 1.24e-13  & 0.0033 & 0.0325 &  0.0235\\
				\hline
	\end{tabular}}
\end{table}

Table~\ref{comp1solkdv} shows that the MC and EC schemes described in Section~\ref{KdVmeth} preserve two conservation laws to machine accuracy. The most accurate of these schemes is $\mbox{EC}_{10}(0.12)$. Minimizing the error in the non-preserved conservation law does not optimise the numerical solution, but nevertheless yields a solution error that is comparable to (for MC$_{10}$) or smaller than (for MC$_{8}$ and EC$_{10}$) the errors in the Furihata, multisymplectic and narrow box schemes.

\begin{figure}[htbp]
	\centering\includegraphics[width=14cm,height=17cm]{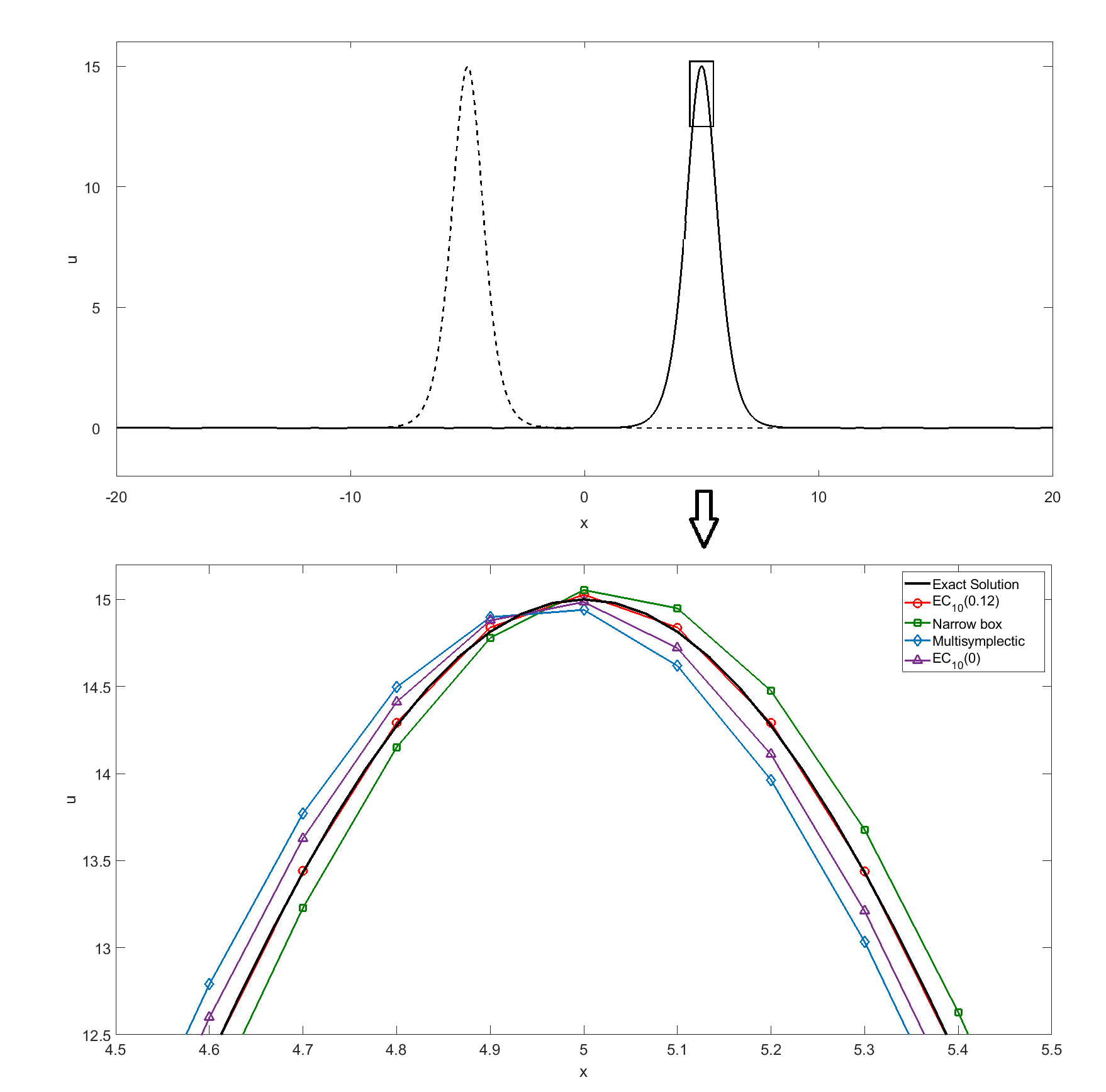}
	\caption{One-soliton solution for the KdV equation (\ref{KdV}) with initial condition (\ref{initsolKdV}) and $c=d=5$; the stepsizes are $\Delta x=0.1$ and $\Delta t=0.01$. Top: 
		Initial condition (dashed curve) and numerical solution $\mbox{EC}_{10}(0.12)$ at $T=2$ (solid curve). Bottom: Top of the soliton at $T=2$; exact solution and numerical solutions from $\mbox{EC}_{10}(0.12)$, $\mbox{EC}_{10}(0)$, multisymplectic and narrow box schemes.}
	\label{1solfigKdVnew}
\end{figure}

The upper part of Fig.~\ref{1solfigKdVnew} shows the initial condition and the numerical solution $\mbox{EC}_{10}(0.12)$ at the final time $T=2$. The lower plot shows only the top of the soliton, comparing the exact solution (\ref{ex1solKdV}) with the numerical solutions given by $\mbox{EC}_{10}(0.12)$, the multisymplectic, narrow box and Furihata ($\mbox{EC}_{10}(0)$) schemes. The $\mbox{EC}_{10}(0.12)$ solution is the closest to the exact solution, reflecting the results in Table~\ref{comp1solkdv}. 

The second benchmark test is the interaction between two solitons. The exact solution on the infinite line is
\begin{equation}\label{KdV2solex}
u_{\text{exact}}(x,t)=\frac{12(c_1-c_2)\left(c_1\cosh^2\xi_2+c_2\sinh^2\xi_1\right)}{\left((\sqrt{c_1}-\sqrt{c_2})\cosh(\xi_1+\xi_2)+(\sqrt{c_1}+\sqrt{c_2})\cosh(\xi_1-\xi_2)\right)^2},
\end{equation}
where
\begin{equation}\label{KdV2solexpar}
\xi_1=\frac{\sqrt{c_1}}2(x+d_1-c_1t),\qquad  \xi_2=\frac{\sqrt{c_2}}2(x+d_2-c_2t).
\end{equation}
Again, we use step sizes $\Delta x=0.1$ and $\Delta t=0.01$ over the spatial domain $[-20,20]$ with periodic boundary conditions, on the temporal interval $[0,2]$. The initial condition is obtained by evaluating (\ref{KdV2solex})--(\ref{KdV2solexpar}) at $t=0$, using the parameters
\begin{equation}\label{par2solKdV}
c_1=10, \quad c_2=5, \quad d_1=12, \quad d_2=10.
\end{equation}

 For this problem, the values $\alpha_1=-0.099$, $(\alpha_2,\beta_2)=(-0.011,-0.031)$ and $\alpha_3=0.23$ minimize the solution error for $\mbox{MC}_8(\alpha_1)$, $\mbox{MC}_{10}(\alpha_2,\beta_2)$ and $\mbox{EC}_{10}(\alpha_3)$, respectively. The values $\alpha_1=-0.084$ and $\alpha_3=0.26$, which minimize the error in the non-preserved conservation law, both produce fairly accurate solutions. The error in the energy conservation law for $\mbox{MC}_{10}(\alpha_2,\beta_2)$ is minimized (but remains large) for each tested $\beta_2$ by a large negative value of $\alpha_2$; this produces a large solution error. So minimising the remaining conservation law is a poor criterion for selecting between the schemes $\mbox{MC}_{10}(\alpha_2,\beta_2)$.

\begin{table}[t!]
\caption{Errors in conservation laws, solution and phase shift at the final time $T=2$ for the two-soliton problem with parameters {\rm(\ref{par2solKdV})} over $[-20,20]$ with periodic boundary conditions using different methods with $\Delta x=0.1$, $\Delta t=0.01$. The non-zero values of the free parameter minimize the error indicated by a star.}
\label{comp2solkdv}
	\small
		\centerline{\begin{tabular}{cccccc}
		\hline
Method &  $\text{Err}_1$ & $\text{Err}_2$ & $\text{Err}_3$  &  Solution error &  $\mbox{Err}_\phi$ \\
\hline
$\mbox{EC}_8$	& 2.27e-13   & 0.0201 & 5.64e-11 &  0.4561 & 0.36\\
$\mbox{MC}_8(0)$	& 3.69e-13 & 5.17e-12 & 30.3753 & 0.3338 & 0.26\\
$\mbox{MC}_8(-0.099)$	& 2.56e-13  & 3.87e-12 & 4.9224 & \phantom{$^*$}0.0301$^*$& -0.04 \\
$\mbox{MC}_8(-0.084)$	& 2.70e-13  & 5.00e-12 & \phantom{$^*$}0.3785$^*$ & 0.0625& 0.06 \\
Furihata; $\mbox{EC}_{10}(0)$	& 1.71e-13  & 1.3595 & 2.18e-11 &  0.1706 & 0.16\\
$\mbox{EC}_{10}(0.23)$	& 1.99e-13  & 0.1725 & 2.00e-11 & \phantom{$^*$}0.0213$^*$ & -0.04 \\
$\mbox{EC}_{10}(0.26)$	& 1.85e-13 & \phantom{$^*$}0.0187$^*$ & 2.55e-11 & 0.0301& -0.04 \\
$\mbox{MC}_{10}(0,0)$	& 1.42e-13  & 1.36e-12 & 27.7429 &  0.2391& 0.16\\
$\mbox{MC}_{10}(-0.011,-0.031)$	& 1.85e-13  & 1.71e-12 & 28.6356 & \phantom{$^*$}0.0253$^*$& -0.04 \\
Multisymplectic & 1.85e-13  & 0.4373 & 28.1328 &  0.2557 & 0.26\\
Narrow box & 1.71e-13  & 0.8633 & 17.9549 &  0.0255& 0.06\\
			\hline
	\end{tabular}}
\end{table}

Table~\ref{comp2solkdv} shows the solution error (\ref{errkdvsol}) and the error in the three conservation laws according to (\ref{errkdvclCLP}) or, for non-preserved conservation laws, (\ref{errkdvNCLP}). The table includes the error in the phase shift for the fastest soliton at the final time,
\begin{equation*}\mbox{Err}_\phi=(x_{\text{max}}-\tilde{x}_{\text{max}}){\vert_{t=2}},
\end{equation*}
where $x_{\text{max}}$ and $\tilde{x}_{\text{max}}$ denote the location of the peak of the fastest soliton in the exact and numerical solution respectively. 

Table~\ref{comp2solkdv} shows that $\mbox{MC}_{8}(-0.099)$, $\mbox{EC}_{10}(0.23)$, $\mbox{MC}_{10}(-0.011,-0.031)$ and the narrow box scheme give the most accurate solutions. The schemes obtained by choosing the free parameter in $\mbox{MC}_8$ and $\mbox{EC}_{10}$ to minimize the error in the non-preserved conservation law are more accurate than the Furihata and multisymplectic schemes.

\begin{figure}[htbp]
\centering{	\includegraphics[width=14cm,height=17cm]{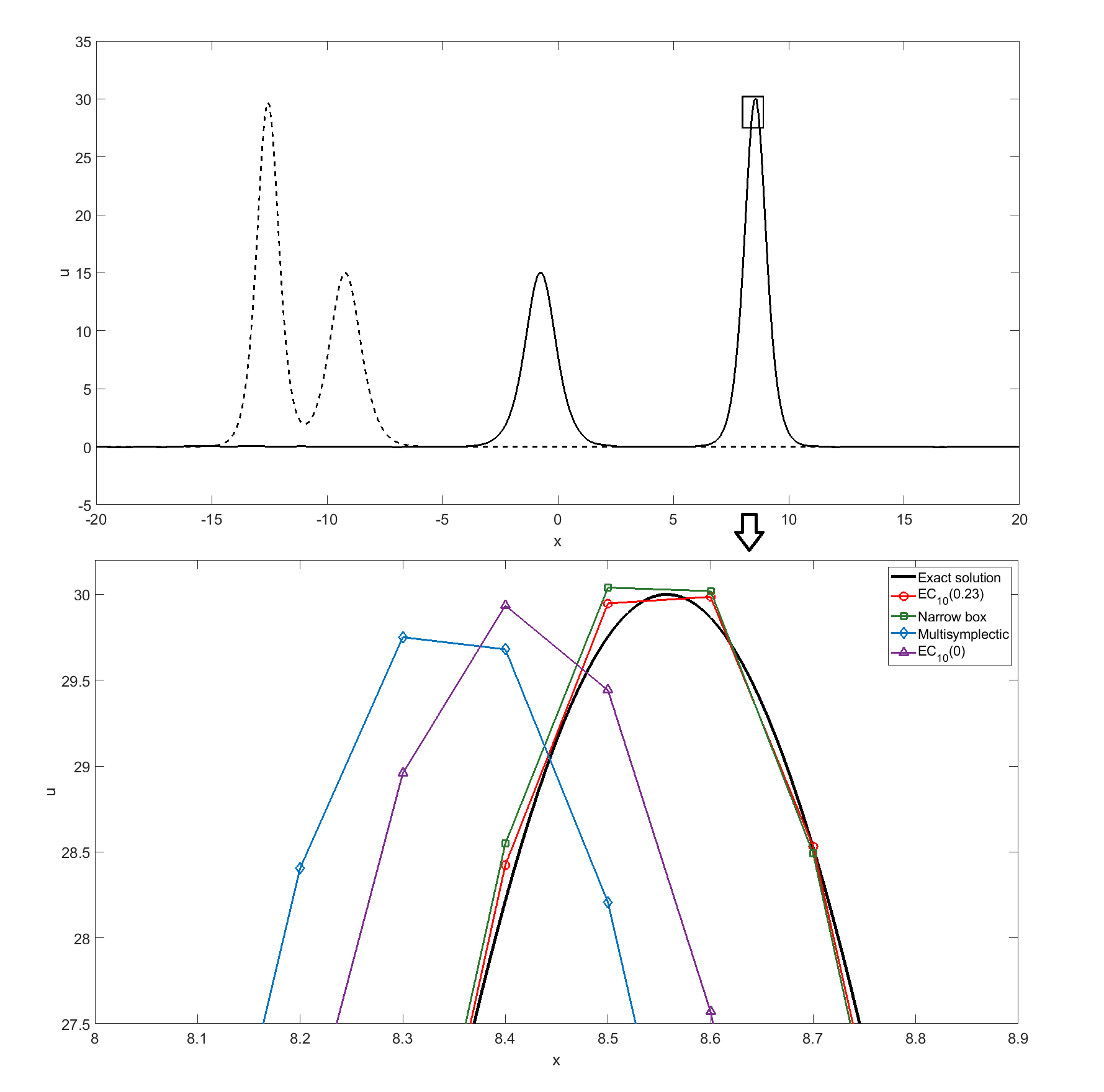}}
	\caption{Two-soliton solutions of KdV with the parameters (\ref{par2solKdV}) and stepsizes $\Delta x=0.1$ and $\Delta t=0.01$. Top: 
		Initial condition (dashed curve) and numerical solution given by $\mbox{EC}_{10}(0.23)$ at $T=2$ (solid curve). Bottom: The top of the fastest soliton at $T=2$; exact solution and numerical solutions from the $\mbox{EC}_{10}(0.23)$, $\mbox{EC}_{10}(0)$, multisymplectic and narrow box schemes.}
	\label{ft2KdV}
\end{figure}

The upper part of Figure~\ref{ft2KdV} shows the initial condition (dashed line) and the numerical solution $\mbox{EC}_{10}(0.23)$ at time $T=2$ (solid line). The lower plot shows the exact solution (\ref{KdV2solex}) and the numerical solutions from $\mbox{EC}_{10}(0.23)$, the Furihata ($\mbox{EC}_{10}(0)$), multisymplectic and narrow box schemes. The narrow box and $\mbox{EC}_{10}(0.23)$ schemes are the most accurate and give similar results.

As a last numerical test, we solve the two-soliton problem on a coarser time grid, setting $\Delta x=0.1$ and $\Delta t=0.02$. 
For these stepsizes, the values $\alpha_1=-0.22$, $(\alpha_2,\beta_2)=(0.815,-0.001)$ and $\alpha_3=0.66$ minimise the solution error for $\mbox{MC}_8(\alpha_1)$, $\mbox{MC}_{10}(\alpha_2,\beta_2)$ and $\mbox{EC}_{10}(\alpha_3)$. The values $\alpha_1=-0.16$
and $\alpha_3=0.53$ yield the minimal error in the non-preserved conservation law. 
\begin{table}[t!]
	\caption{Errors in conservation laws, solution and phase shift at the final time $T=2$, for the two-soliton KdV problem with parameters {\rm(\ref{par2solKdV})} over $[-20,20]$ with periodic boundary conditions, using different methods with $\Delta x=0.1$, $\Delta t=0.02$. The non-zero values of the free parameter minimize the error indicated by a star.}
	\small
		\centerline{\begin{tabular}{cccccc}
		\hline
			Method &  $\text{Err}_1$ & $\text{Err}_2$ & $\text{Err}_3$  &  Solution error &  $\text{Err}_\phi$\\
			\hline
$\mbox{EC}_8$	& 2.13e-13   & 1.2732 & 6.37e-11 &  0.6704& 0.66\\
$\mbox{MC}_8(0)$	& 1.99e-13 & 1.99e-12 & 55.4259 & 0.6281& 0.56\\
$\mbox{MC}_8(-0.22)$	& 3.98e-13  & 5.97e-12 & 15.8171 & \phantom{$^*$}0.1065$^*$& 0.06 \\
$\mbox{MC}_8(-0.16)$	& 2.27e-13  & 4.21e-12 & \phantom{$^*$}5.997$^*$ & 0.2176& 0.16 \\
Furihata; $\mbox{EC}_{10}(0)$	& 1.99e-13  & 2.7234 & 2.82e-11 &  0.4501& 0.36\\
$\mbox{EC}_{10}(0.66)$	& 1.42e-13  & 0.5161 & 2.00e-11 & \phantom{$^*$}0.0749$^*$ & 0.06 \\
$\mbox{EC}_{10}(0.53)$	& 1.14e-13 & \phantom{$^*$}0.1326$^*$ & 2.73e-11 & 0.1245& 0.06 \\
$\mbox{MC}_{10}(0,0)$	& 1.56e-13  & 1.82e-12 & 54.2290 &  0.5621 & 0.46\\
$\mbox{MC}_{10}(0.815,-0.001)$	& 1.56e-13  & 1.25e-12 & 59.5767 & \phantom{$^*$}0.0947$^*$ & 0.06\\
Multisymplectic & 2.27e-13  & 0.4306 & 53.2081 &  0.5678 & 0.46\\
Narrow box & 2.13e-13  & 0.8481 & 10.3316 &  0.3860& 0.36\\
			\hline
	\end{tabular}}
	\label{comp2solkdv2}
\end{table}
\begin{figure}[h!]
	\centering{\includegraphics[width=15cm,height=8cm]{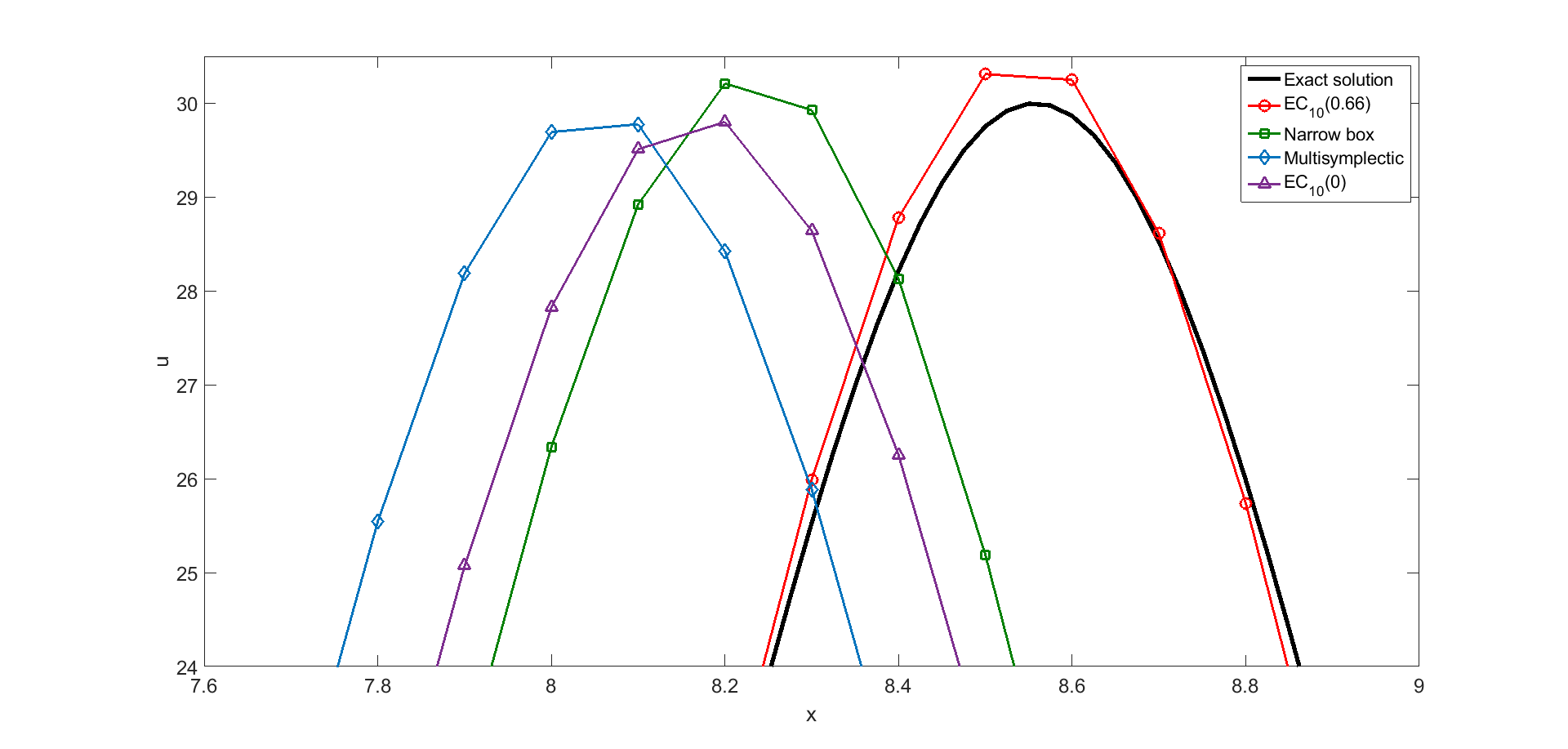}}
	\caption{Two-soliton solutions for (\ref{KdV}) with the parameters (\ref{par2solKdV}) and stepsizes $\Delta x=0.1$ and $\Delta t=0.02$. Top of the fastest soliton at time $T=2$: Exact solution and numerical solutions from the $\mbox{EC}_{10}(0.66)$, $\mbox{EC}_{10}(0)$, multisymplectic and narrow box schemes.}
	\label{ft2KdV2}
\end{figure}

Most results in Table \ref{comp2solkdv2} are qualitatively similar to their counterparts in Table \ref{comp2solkdv}, though with larger solution and phase errors. However, the narrow box scheme is far less accurate for the larger time step. The solution error in the most accurate scheme, $\text{EC}_{10}(0.66)$, is around $3.5$ times that of the most accurate $\text{EC}_{10}$ scheme for the smaller time step $\Delta t=0.01$. This growth in solution error is slightly greater that those of the Furihata ($\text{EC}_{10}(0)$) and multisymplectic schemes (whose phase errors also grow more slowly). Even so, $\text{EC}_{10}(0.66)$ is by far the most accurate of the schemes (see Fig. \ref{ft2KdV2}).

\section{A nonlinear heat equation}\label{heat}

In this section we consider the nonlinear heat equation,
\begin{equation}\label{heateq}
\mca\equiv u_t-u_x^2-uu_{xx}=0,\qquad (x,t)\in\Omega~\equiv~[a,b]\times[0,\infty),
\end{equation}
coupled with suitable initial and boundary conditions:
\begin{equation}\label{initbound}
u(x,0)=\psi(x),\qquad u(a,t)=\varphi_1(t),\qquad u(b,t)=\varphi_2(t).
\end{equation}
Equation (\ref{heateq}) has only two independent (equivalence classes of) conservation laws: 
\begin{eqnarray}\label{CL1heat}
\mca &=& D_t(G_1)+D_x(F_1)\equiv D_t(u) + D_x(-uu_x),\\ \label{CL2heat}
x\mca &=& D_t(G_2)+D_x(F_2)\equiv D_t(xu) + D_x\left({u^2}/{2}-xuu_x\right),
\end{eqnarray}
with characteristics 
\begin{equation}\label{heatchar}
\mcq_1=1,\quad\mcq_2=x,
\end{equation}
respectively \citep[see][]{CRC}. To construct finite difference schemes that preserve a discrete version of (\ref{CL1heat}) and (\ref{CL2heat}), we use the following general results.

\begin{theo}\label{theoheat}
	
Any partial differential equation of the form 
\begin{equation}\label{genPDEth}
D_t\left(g[u]\right)+D_x^{\,k}\left(f[u]\right)=0,\qquad k\in\mathbb{N},
\end{equation}
where $g[u]$ and $f[u]$ are smooth functions of $u$ and its derivatives, has $k$ conservation laws whose characteristics are $Q_i=x^{i-1},\, i=1,\ldots,k.$ Any scheme of the form
\begin{equation}\label{genschemeheat}
D_n\left(\widetilde{g[u]}\right)+D_m^{\,k}\left(\widetilde{f[u]}\right)=0,
\end{equation}
where $\widetilde{g[u]}$ and $\widetilde{f[u]}$ are finite difference approximations to $g[u]$ and $f[u]$ respectively, has $k$ conservation laws whose characteristics are $\widetilde{Q_i}=x_0^{i-1},\, i=1,\ldots,k$.\\\quad

\begin{proof} On solutions of (\ref{genPDEth}), using integration by parts,
$$x^{i-1}\left( D_t(g[u])+D_x^{\,k}(f[u])\right)=D_t\left(x^{i-1}g[u]\right)+(-1)^kD_x^{\,k}(x^{i-1})f[u]+D_x(h[u]),$$
for some function $h$ of $u$ and its derivatives. As $i-1<k$, this simplifies to
$$x^{i-1}\left( D_t(g[u])+D_x^{\,k}(f[u])\right)=D_t\left(x^{i-1}g[u]\right)+D_x(h[u]),$$
which is a divergence and, therefore, a conservation law with characteristic $Q_i=x^{i-1}.$ 

The proof for the discrete case is similar. Summation by parts gives
$$x_0^{i-1}\left( D_n\left(\widetilde{g[u]}\right)+D_m^{\,k}\left(\widetilde{f[u]}\right)\right)=D_n\left(x_0^{i-1}\widetilde{g[u]}\right)+\left(-S_m^{-1}D_m\right)^k(x_0^{i-1})\widetilde{f[u]}+D_m\left(\widetilde{h[u]}\right)$$
on solutions of (\ref{genschemeheat}), for some function $\widetilde{h[u]}$ of $u$ and its shifts. Again $i-1<k$, so
$$x_0^{i-1}\left( D_n\left(\widetilde{g[u]}\right)+D_m^{\,k}\left(\widetilde{f[u]}\right)\right)=D_n\left(x_0^{i-1}\widetilde{g[u]}\right)+D_m\left(\widetilde{h[u]}\right),$$
which is a discrete conservation law with characteristic $\widetilde{Q_i}=x_0^{i-1}.$

\end{proof}
\end{theo}

\subsection{Conservative methods for the nonlinear heat equation}\label{heatmeth}

\begin{figure}
\center{
\begin{tikzpicture}
\draw[fill] (2,0) circle [radius=0.075];
\draw[fill] (4,0) circle [radius=0.075];
\draw[fill] (6,0) circle [radius=0.075];
\draw[fill] (2,2) circle [radius=0.075];
\draw[fill] (4,2) circle [radius=0.075];
\draw[fill] (6,2) circle [radius=0.075];
\node [below right] at (4,1) {\tiny{$(0,1/2)$}};
\node [below] at (3,1) {\tiny{$(-1/2,1/2)$}};
\draw (4,1) circle [radius=0.175];
\draw (2.87,0.87)--(3.13,1.13);
\draw (2.87,1.13)--(3.13,0.87);
\draw (3.87,-0.13)--(4.13,0.13);
\draw (3.87,0.13)--(4.13,-0.13);
\draw (3.87,-0.13) rectangle (4.13,0.13);
\node [below right] at (4,0) {\tiny{$(0,0)$}};

\end{tikzpicture}
\caption{The 6-point rectangular stencil for (\ref{heateq}). Conservation laws are preserved to second-order at the central point $(0, 1/2)$ (circle), densities and fluxes respectively at $(0,0) $ and $(-1/2, 1/2)$ (crosses).}
\label{stencilheat}
}
\end{figure}
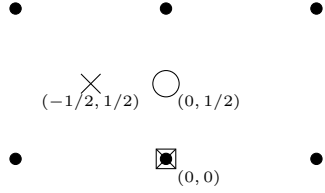

The nonlinear heat equation (\ref{heateq}) is of the form (\ref{genPDEth}), with $k=2$. Therefore, according to Theorem~\ref{theoheat}, both conservation laws can be preserved by finding suitable finite difference approximations of $g[u]=u$ and $f[u]=-{u^2}/{2}$. This can be achieved to second order on the most compact rectangular stencil for (\ref{heateq}), which consists of six points (see Fig.~\ref{stencilheat}).

This results in a two-parameter family of schemes,
\begin{equation}\label{CSab}
\mbox{CS}(\alpha,\beta)\equiv\widetilde{\mca}\equiv D_m\widetilde{F_1}+D_n\widetilde{G_1}=0,
\end{equation}
with
\begin{align}\label{fluxheat}
\widetilde{F_1}&=D_m\left(\widetilde{f[u]}\right),\quad\text{where}\quad\widetilde{f[u]}=-\frac{1}2u_{-1,0}u_{-1,1}+\Delta t^2\beta(D_nu_{-1,0})^2,\\\label{densheat}\widetilde{G_1}&=\widetilde{g[u]}=u_{0, 0}+\alpha \Delta x^2D_m^2u_{-1, 0}\,.
\end{align}
These schemes preserve a discrete version of (\ref{CL1heat}) with $\widetilde{Q_1}=1$ and of (\ref{CL2heat}), namely
\begin{equation}\label{CSabCL}
\widetilde{Q_2}\widetilde\mca=D_m\widetilde{F_2}+D_n\widetilde{G_2},
\end{equation}
where
\[
\widetilde{Q_2}=x_0\,,\qquad
\widetilde{F_2}=\mu_m(x_{-1})\widetilde{F_1}-\mu_m\left(\widetilde{f[u]}\right),\qquad
\widetilde{G_2}=x_0\widetilde{G_1}\,.
\]
Except for $\mbox{CS}(0,0)$, the approximated densities and fluxes include derivative terms that do not appear in the corresponding continuous quantities. However, these vanish as the stepsizes approach zero. The following schemes are particularly straightforward.

The scheme $\mbox{CS}(0,-1/4)$ has perhaps the most obvious discretization of $\widetilde{F_1}$ and $\widetilde{G_1}$, namely
\[
\widetilde{F_1}=D_m\left(-\tfrac{1}{2}\mu_nu_{-1,0}^2\right),\qquad
\widetilde{G_1}=u_{0, 0}\,.
\]
The components of the second conservation law are
\[
\widetilde{F_2}=\mu_m(x_{-1})\widetilde{F_1}+\tfrac{1}2\mu_m\left(\mu_nu_{-1,0}^2\right),\qquad
\widetilde{G_2}=x_0u_{0,0}\,.
\]
The scheme $\mbox{CS}(0,-1/8)$ amounts to
\begin{equation*}
D_n(u_{0,0})+D_m^2\left(-\tfrac{1}{2}\left(\mu_nu_{-1,0}\right)^2\right)=0,
\end{equation*}
which is obtained by applying the implicit midpoint rule to the following simple semidiscretization of (\ref{heateq}), with $U_i(t)$ approximating $u(x_i,t)$:
$$U_0'+D_m^2\left(-\tfrac{1}2(U_{-1})^2\right)=0.$$
The densities and fluxes of the discrete conservation laws preserved by $\mbox{CS}(0,-1/8)$ are
\[
\widetilde{F_1}=D_m\left(-\tfrac{1}{2} \left(\mu_nu_{-1,0}\right)^2\right),\qquad
\widetilde{G_1}=u_{0, 0}\,;
\]
\[
\widetilde{F_2}=\mom(x_{-1})\widetilde{F_1}+\tfrac{1}2\mom\left(\left(\mu_nu_{-1,0}\right)^2\right),\qquad
\widetilde{G_2}=x_0u_{0,0}\,.
\]
The scheme with the fewest terms is $\mbox{CS}(0,0)$, whose densities and fluxes are 
\[
\widetilde{F_1}=D_m\left(-\tfrac{1}{2}{u_{-1,0}u_{-1,1}}\right),\qquad
\widetilde{G_1}=u_{0, 0}\,;
\]
\[
\widetilde{F_2}=\mu_m(x_{-1})\widetilde {F_1}+\tfrac{1}2\mu_m\left({u_{-1,0}u_{-1,1}}\right),\qquad
\widetilde{G_2}=x_0u_{0,0}\,.
\]
This scheme can be solved explicitly at the $k^{\mathrm{th}}$ integration step, provided that the matrix
\begin{equation*}
\frac{1}{\Delta t}\,I-\frac{1}{2\Delta x^2}\,M\circ V_{k-1}
\end{equation*}
is invertible; here $I$ is the identity matrix, $\circ$ is the Hadamard product,\footnote{The Hadamard product gives $(A\circ B)_{i,j}=A_{i,j}B_{i,j},$ for each $i,j.$}
$$M=\left(\begin{array}{ccccc}
-2& 1 & & &\\
1 & -2 & 1 & &\\
& 1 & -2 & 1 & \\
&  & \ddots & \ddots & \ddots 
\end{array}\right),\qquad V_{k-1}=\mathbf{e}\,u_{k-1}\,,$$ where $\mathbf{e}$ is the column vector whose entries are all $1$, and $u_{k-1}$ is the row vector whose entries are the approximation (from the previous step) of $u$ at the spatial grid points.

\begin{rem}
The nonlinear heat equation (\ref{heateq}) is a special case of the porous medium equation, 
\begin{equation}\label{PME}
D_t(u)+D_x^2\left(-\,\frac{u^s}{s}\right)=0,\qquad (x,t)\in\Omega~\equiv~[a,b]\times[0,\infty),
\end{equation}
where $s\in \mathbb{N}\backslash\{1\}$. For each $s$, this equation has only two conservation laws, with characteristics $1$ and $x$. The approach that we have used for the nonlinear heat equation can be used to obtain conservative schemes for (\ref{PME}) with $s>2$.  On the six-point stencil in Fig.~\ref{stencilheat}, this yields an $s$-parameter family of second-order methods.
\end{rem}
\subsection{Numerical tests}

In this section, three benchmark numerical tests for the problem (\ref{heateq})-(\ref{initbound}) are used to show the effectiveness of the methods developed in Section~\ref{heatmeth}. We compare the results from several $\mbox{CS}(\alpha,\beta)$ schemes, which preserve both conservation laws, with those from the following second-order scheme that, in general, does not preserve either conservation law: 
\begin{align}\label{stand2}
D_nu_{0,0}-(D_m\mu_m\mu_nu_{-1,0})^2-\mu_n(u_{0,0})D_m^2(\mu_nu_{-1,0})=0.
\end{align}
We call the scheme ML/IM, as it is obtained by applying the implicit midpoint method to the following standard second-order method-of-lines semidiscretization of (\ref{heateq}):
\begin{align*}
U_0'-(D_m\mu_mU_{-1})^2-U_0D_m^2U_{-1}=0.
\end{align*}

Again, implicit methods are solved by a simplified Newton method with frozen Jacobian, run until the error reaches machine accuracy. The (relative) solution error at the final time $t=T$ is evaluated as
\begin{equation}\label{relerr}
\frac{\|u-u_{\text{exact}}\|}{\|u_{\text{exact}}\|}\bigg\vert_{t=T}.
\end{equation}
The errors in the discrete conservation laws (\ref{CSab}) and (\ref{CSabCL}) are evaluated respectively as\footnote{If periodic or zero boundary conditions apply, (\ref{errCL1}) and (\ref{errCL2}) can be replaced with
	$${\Delta x}\max_{j=1,\ldots, N-1} \sum_{i=1}^{M}D_n(u_{i,1}+\alpha\Delta x^2D_m^2u_{i-1,j}),\qquad{\Delta x}\max_{j=1,\ldots, N-1} \sum_{i=1}^{M}{x_i}D_n\left(u_{i,j}+\alpha\Delta x^2D_m^2u_{i-1,j}\right),$$
	which measure the error in the conservation of the global invariants $$\int G_1\,\mathrm d x = \int u \,\mathrm d x \quad \mbox{and} \quad \int G_2\,\mathrm d x = \int xu \,\mathrm d x.$$}
\begin{align}\label{errCL1}
\mbox{Err}_1=\,\Delta x\!\max_{j=1,\ldots, N-1} &\left|\,\sum_{i=2}^{M-1}D_n\!\left(u_{i,j}+\alpha\Delta x^2D_m^2u_{i-1,j}\right)\!-\!D_m\left\{\left[\frac{u_{r,j}u_{r,j+1}+2\beta\Delta t^2(D_nu_{r,j})^2}{2\Delta x}\right]_{r=1}^{M-1}\right\}\right|\,,\\\nonumber
\mbox{Err}_2=\,\Delta x\!\max_{j=1,\ldots, N-1} &\left|\,\sum_{i=2}^{M-1}x_iD_n\!\left(u_{i,j}+\alpha\Delta x^2D_m^2u_{i-1,j}\right)\!-\!\left[\frac{\mu_m(x_r)D_m(u_{r,j}u_{r,j+1})-\mu_m(u_{r,j}u_{r,j+1})}{2\Delta x}\right]_{r=1}^{M-1}\right.\\
&\left.+\beta\Delta t^2\left[\frac{\mu_m(x_r)D_m\!\left((D_nu_{r,j})^2\right)-\mu_m\left((D_nu_{r,j})^2\right)}{\Delta x}\right]_{r=1}^{M-1}\,\,
\right|\,,\label{errCL2}
\end{align}
where $u_{i,j}\simeq u(a+i\Delta x, j\Delta t)$, so that subscripts denote shifts with respect to $(x,t)=(a,0)$. To evaluate the error in the conservation laws resulting from ML/IM we use (\ref{errCL1}) and (\ref{errCL2}), setting $\alpha=\beta=0$.

The first benchmark problem is (\ref{heateq}) with the initial and boundary conditions
\begin{align}\label{init3}
u(x,0)=\left(1-\frac{x^2}{6}\right)_{\!\!+}, \quad -6\leq x\leq 6,\qquad u(-6,t)=u(6,t)=0,\quad t\in[0,4],
\end{align}
where $f_+=\max(f,0)$. These conditions yield the Barenblatt solution of the porous medium equation (\ref{PME}) with $s=2$, which is
\begin{equation*}
u_{\text{exact}}(x,t)=(t+1)^{-1/3}\left(1-\frac{x^2}{6(t+1)^{2/3}}\right)_{\!\!+}.
\end{equation*}
For all $t>0$, this solution has compact support with the interface moving outward at a finite speed. The Barenblatt solution is a (weak) energy solution, but not a classical solution as it is not differentiable at the interface points. Such solutions cause difficulties in numerical simulation. Standard finite element methods can create oscillations close to the interface, but negative values have no meaning physically \citep[see][]{ZW08}. Here we show that, by contrast, various conservative finite difference schemes $\text{CS}(\alpha,\beta)$ are effective for non-smooth solutions. For simplicity, we will consider only the one-parameter family obtained by setting $\alpha=0$.

Table~\ref{errBern1} shows the errors in the conservation laws for various $\text{CS}(0,\beta)$ given the stepsizes $\Delta x=0.25$ and $\Delta t=0.333$. These schemes locally preserve both conservation laws to machine accuracy. The solution error at the final time $T=4$, evaluated according to (\ref{relerr}), is minimised by setting $\beta=0.21$. Nevertheless, the explicitly-solved scheme $\mbox{CS}(0,0)$, though slightly less accurate, is a better option because of its low computation time.

ML/IM does not converge on such a coarse grid. Only by reducing the time step so that $\Delta t<\Delta x^2$ can this scheme be made to converge. Reducing the time step to $\Delta t=0.0267$, the solution error is still larger than those of the $\mbox{CS}(0,\beta)$ methods, which converge even when $\Delta t>\Delta x$. Note that ML/IM preserves the conservation law (\ref{CL2heat}); this is a consequence of the reflectional symmetry of the scheme and the boundary conditions.

Tables~\ref{errBern2} and \ref{errBern3} show the outcomes of solving the same problem with various $\mbox{CS}(0,\beta)$ on the finer grids $\Delta x=0.1$, $\Delta t=0.133$ and $\Delta x=0.025$, $\Delta t=0.03$. On these grids the values $\beta=0.07$ and $\beta=0.05$ respectively minimize the solution error. The explicit scheme $\mbox{CS}(0,0)$ is by far the most efficient and has a low solution error. Of the implicit schemes, the optimised scheme is the fastest in each case. The errors in the conservation laws are tiny, but grow as the grid is refined due to the accumulation of rounding errors.
Again, ML/IM requires smaller timesteps for convergence; even then, the solution error is still far greater than those of the conservative methods.

\begin{table}[t!]
\caption{Errors in solution and conservation laws at the final time $T=4$, when solving {\rm(\ref{heateq})} with the conditions {\rm(\ref{init3})}, using {\rm\mbox{CS}}$(0,\beta)$ schemes with $\Delta x=0.25$, $\Delta t=0.333$ and the scheme {\rm ML/IM} with $\Delta x=0.25$, $\Delta t=0.0267$.}
\small
	\centerline{\begin{tabular}{ccccc}
\hline
Method &  $\text{Err}_1$ & $\text{Err}_2$ & Solution error & Computational time\\ 
\hline
$\mbox{CS}(0,-1/4)$	& 8.11e-16   & 5.66e-16 & 0.0038 &  0.032 \\
$\mbox{CS}(0,-1/8)$	&  4.55e-16   & 3.92e-16 & 0.0035 & 0.026\\
$\mbox{CS}(0,0)$	&  5.95e-15   & 5.41e-15 & 0.0032 & 0.002 \\
$\mbox{CS}(0,0.21)$	&  8.61e-16   & 1.03e-15 & 0.0028 & 0.030 \\
ML/IM &  0.0671   & 5.30e-15 & 0.0307 & 0.221\\
\hline
\end{tabular}}
\label{errBern1}
\end{table}
\begin{table}[t!]
\caption{Errors in solution and conservation laws at the final time $T=4$, when solving {\rm(\ref{heateq})} with the conditions in {\rm(\ref{init3})}, using {\rm CS}$(0,\beta)$ methods with $\Delta x=0.1$, $\Delta t=0.133$ and the scheme {\rm ML/IM} with $\Delta x=0.1$, $\Delta t=0.005$.}
\small
	\centerline{\begin{tabular}{ccccc}
	\hline
Method &  $\text{Err}_1$ & $\text{Err}_2$ &  Solution error & Computational time\\
\hline 
$\mbox{CS}(0,-1/4)$	& 1.66e-15   & 2.57e-15 & 0.0013 &  0.15  \\
$\mbox{CS}(0,-1/8)$	&  1.58e-15   & 1.83e-15 & 0.0012 & 0.10\\
$\mbox{CS}(0,0)$	&  4.88e-14   & 4.67e-14 & 0.0011 & 0.014 \\
$\mbox{CS}(0,0.07)$	&  3.70e-15   & 3.36e-15 & 9.77e-04 & 0.09 \\
ML/IM &  0.0232   & 2.37e-14 & 0.0126 & 6.48\\ 
\hline
\end{tabular}}
\label{errBern2}
\end{table}
\begin{table}[t!]
\caption{Errors in solution and conservation laws at the final time $T=4$, when solving {\rm(\ref{heateq})} with the conditions in {\rm(\ref{init3})}, using {\rm \mbox{CS}}$(0,\beta)$ methods with $\Delta x=0.025$, $\Delta t=0.03$ and the scheme {\rm ML/IM} with $\Delta x=0.025$, $\Delta t=0.000267$.}
\small
	\centerline{\begin{tabular}{ccccc}
\hline
Method &  $\text{Err}_1$ & $\text{Err}_2$ &  Solution error & Computational time\\ 
\hline
$\mbox{CS}(0,-1/4)$	& 3.16e-14   & 1.78e-14 & 6.51e-05 &  14.35  \\
$\mbox{CS}(0,-1/8)$	&  2.18e-14   & 1.50e-14 & 5.87e-05 & 8.99\\
$\mbox{CS}(0,0)$	&  3.22e-14   & 2.69e-14 & 5.48e-05 & 0.38 \\
$\mbox{CS}(0,0.05)$	&  3.20e-14   & 4.49e-14 & 5.42e-05 & 6.71 \\
ML/IM &  0.0095   & 2.50e-13 & 0.0035 & 1612.61\\
\hline
\end{tabular}}
\label{errBern3}
\end{table}

\begin{figure}[htbp]
\centering{\includegraphics[width=14cm,height=17cm]{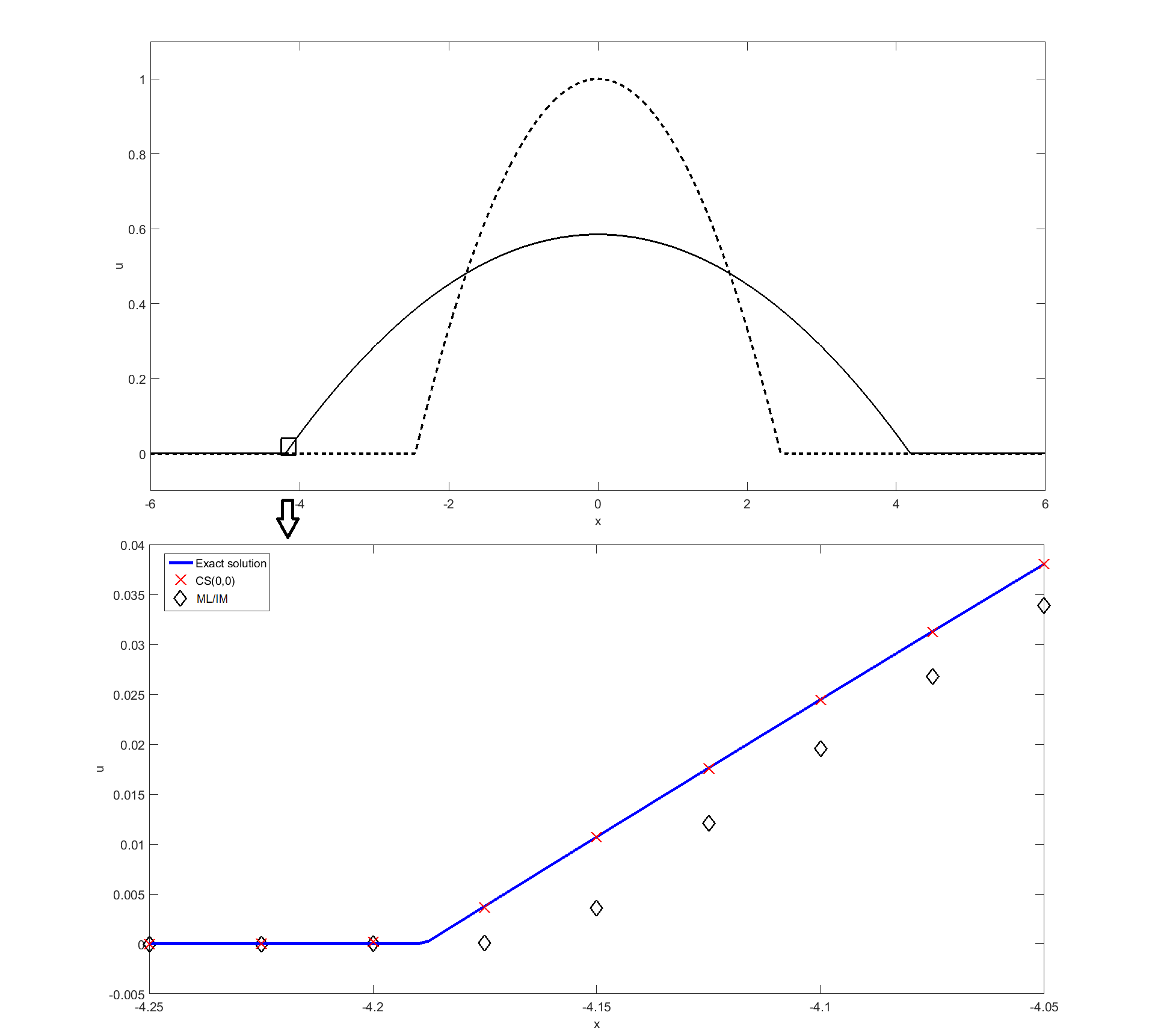}}
\caption{ Results for (\ref{heateq}) with conditions in (\ref{init3}). Top: Initial condition (dashed curve) and numerical solution from $\mbox{CS}(0,0)$ (solid curve), setting $\Delta x=0.025$ and $\Delta t=0.03$ over $[-6,6]$ at time $T=4$. Bottom: Left interface: exact solution (solid curve), numerical solutions $\mbox{CS}(0,0)$ with $\Delta x=0.025$ and $\Delta t=0.03$ (crosses), and ML/IM with $\Delta x=0.025$ and $\Delta t=0.000267$ (diamonds).}
\label{fig3}
\end{figure}

The upper part of Fig.~\ref{fig3} shows the initial condition (dashed line) and the numerical solution given by $\mbox{CS}(0,0)$ for equation (\ref{heateq}), with conditions in (\ref{init3}), setting $\Delta x=0.025$ and $\Delta t=0.03$. The method does not produce any spurious oscillations close to the interface. Magnifying the left interface, as shown at bottom of Fig.~\ref{fig3}, one can see that the solution of $\mbox{CS}(0,0)$ is closer to the exact solution at time $T=4$ than the solution given by ML/IM, even though the time step used to advance ML/IM is much smaller. Furthermore, the interface of the numerical solution has moved at the correct speed and overlaps the interface of the exact solution. The solutions given by $\mbox{CS}(0,\beta)$ for the optimal value of $\beta$ overlap the $\mbox{CS}(0,0)$ solution, so we omit the corresponding figures.

The second benchmark problem is (\ref{heateq}) with $(x,t)\in\Omega=[0,15]\times[0,10]$ and the following initial and boundary conditions:
\begin{align}\nonumber
&u(x,0)=0,\quad x\in[0,15],\\\label{init1}
u(0,t)=t,\quad t\in&\,[0,10],\qquad u(15,t)=0,\quad t\in[0,10].
\end{align}
The exact solution of this problem \citep[see][]{CRC} is again not smooth:
\begin{equation*}
u_{\text{exact}}(x,t)=\begin{cases}
t-x,&\quad 0\leq x \leq t,\\
0,&\quad x> t;
\end{cases}
\end{equation*}
this is a wave travelling with unit speed into an undisturbed medium.

Table~\ref{num2tab1} shows the errors for various $\mbox{CS}(0,\beta)$ with $\Delta x=0.375$ and $\Delta t=0.333$. The value $\beta=-0.14$ gives the minimum solution error. As ML/IM does not converge on this grid, we use the finer timestep $\Delta t=0.00667$ for this method; the problem is not reflectionally symmetric, so neither conservation law is preserved. The results are similar to those for the first benchmark problem. Again, the sub-optimal scheme $\mbox{CS}(0,0)$, solved explicitly, is convenient because of its low computation time. 

\begin{table}[t!]
\caption{Errors in solution and conservation laws at the final time $T=10$, when solving {\rm(\ref{heateq})} with conditions in {\rm(\ref{init1})}, using {\rm{CS}}$(0,\beta)$ schemes with $\Delta x=0.375$, $\Delta t=0.333$ and {\rm ML/IM} with $\Delta x=0.375$, $\Delta t=0.00667$.}
\small
	\centerline{\begin{tabular}{ccccc}
\hline
Method &  $\text{Err}_1$ & $\text{Err}_2$ &  Solution error & Computational time\\
\hline 
$\mbox{CS}(0,-1/4)$	& 2.13e-14   & 3.73e-14 & 0.0013 &  0.06 \\
$\mbox{CS}(0,-1/8)$	&  2.40e-14   & 2.13e-14 & 9.26e-04 & 0.04\\
$\mbox{CS}(0,0)$	&  3.60e-14   & 5.86e-14 & 0.0035 & 0.002 \\
$\mbox{CS}(0,-0.14)$	&  2.40e-14   & 2.66e-14 & 9.13e-04 & 0.04 \\
ML/IM &  0.0845   & 0.7938 & 0.0114 & 1.76\\ 
\hline
\end{tabular}}
\label{num2tab1}
\end{table}
\begin{table}[t!]
\caption{Errors in solution and conservation laws at the final time $T=10$, when solving {\rm(\ref{heateq})} with conditions in {\rm(\ref{init1})}, using {\rm CS}$(0,\beta)$ methods with $\Delta x=0.05$, $\Delta t=0.025$ and {\rm ML/IM} with $\Delta x=0.05$, $\Delta t=8${\rm e-05}.}
\small
	\centerline{\begin{tabular}{ccccc}
\hline
Method &  $\text{Err}_1$ & $\text{Err}_2$ &  Solution error & Computational time\\ 
\hline
$\mbox{CS}(0,-1/4)$	& 9.55e-13   & 1.79e-12 & 1.16e-04 &  4.92  \\
$\mbox{CS}(0,-1/8)$	&  7.92e-13   & 1.34e-12 & 9.99e-05 & 3.67\\
$\mbox{CS}(0,0)$	&  9.96e-13   & 2.06e-12 & 8.16e-05 & 0.19 \\
$\mbox{CS}(0,0.34)$	&  2.60e-12   & 5.37e-12 & 2.94e-05 & 6.90 \\
ML/IM &  0.0114   & 0.1131 & 0.0017 & 1400.70\\ 
\hline
\end{tabular}}
\label{num2tab3}
\end{table}

Similar results are obtained in 
Table~\ref{num2tab3} by solving the same problem using $\mbox{CS}(0,\beta)$ schemes on the finer grid
$\Delta x=0.05$, $\Delta t=0.025$, for which
the value $\beta=0.34$ minimizes the solution error. As ML/IM does not converge on these grids, we
reduce the timestep to $\Delta t=8\text{e-}05$.
Again, the conservative methods are more accurate.

\begin{figure}[htbp]
\centering{\includegraphics[width=14cm,height=17cm]{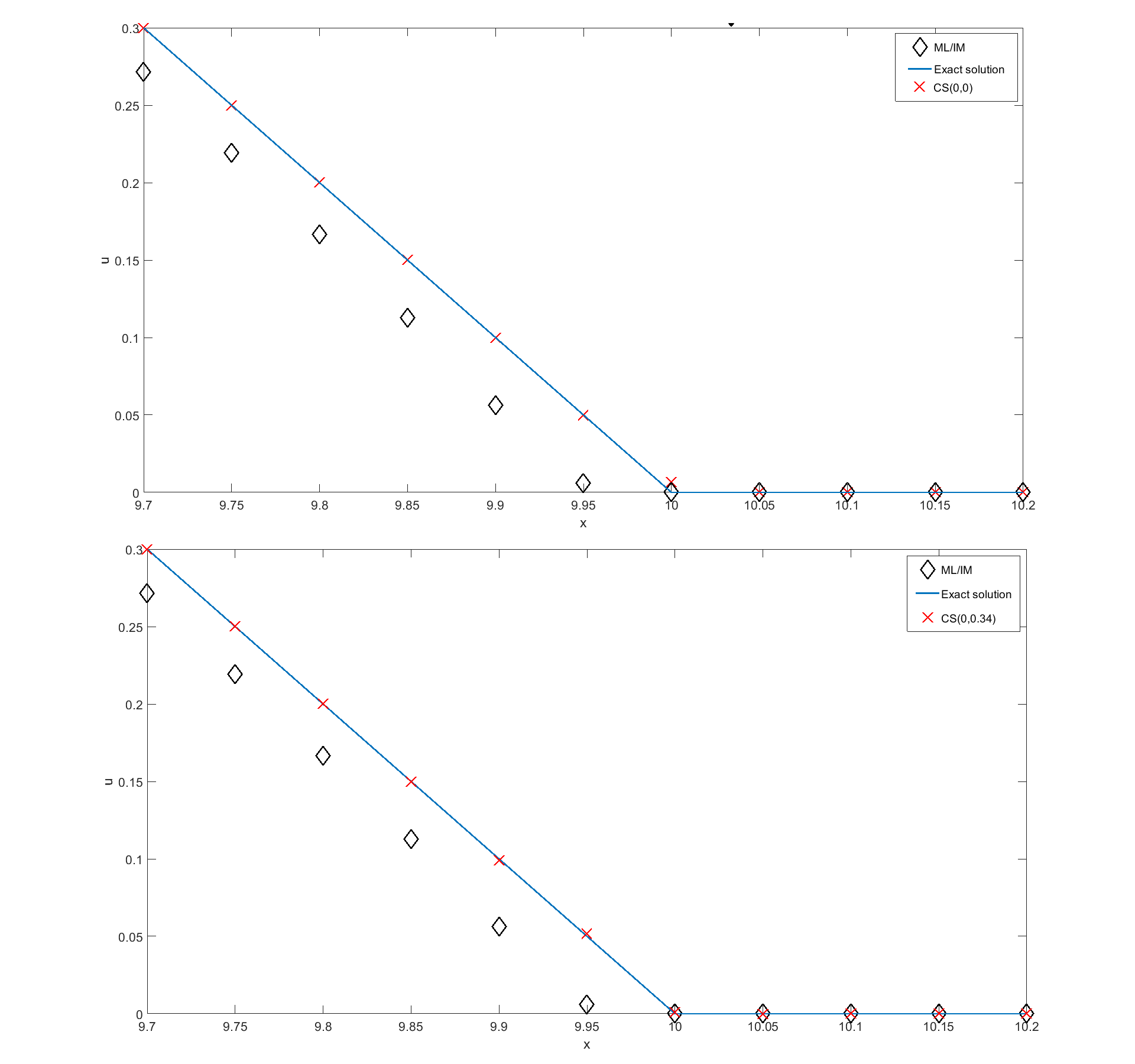}}
\caption{ Results at $T=10$ for equation (\ref{heateq}) subject to (\ref{init1}), setting $\Delta x=0.05$. Exact solution (solid curve), $\text{CS}(0,\beta)$ with $\Delta t=0.025$ (crosses), and ML/IM with $\Delta t=8$e-05 (diamonds). Top: $\beta=0$. Bottom: $\beta=0.34$.}
\label{fig1}
\end{figure}

Figure~\ref{fig1} compares numerical solutions of problem (\ref{heateq}) with (\ref{init1}) on the finer grid, which has a meshpoint at the interface $x=10$. Although ML/IM is qualitatively correct, it produces a slight lag near to the interface. The scheme $\mbox{CS}(0,0)$ is very accurate except at the interface and does not produce spurious oscillations in the solution. The most accurate scheme, $\mbox{CS}(0,0.34)$, models the moving interface extremely well, but produces a very small oscillation in the error at nearby points.

The final benchmark problem is (\ref{heateq}) with the initial and boundary conditions
\begin{align}\nonumber
&u(x,0)=0, \quad 0\leq x\leq 5,\\\label{init2}
&u(0,t)=(\tilde{t}-t)^{-1}\left[1-\left(1-t/\tilde t\,\right)^{2/3}\right],\qquad u(5,t)=0,\quad t\in[0,T]\subset[0,\tilde t\,),
\end{align}
where $\tilde{t}$ is a positive constant (the time of existence of the solution). The exact solution of this problem is \citep[see][]{GP88,CRC}
\begin{equation*}
u_{\text{exact}}(x,t)=\begin{cases}
(\tilde{t}-t)^{-1}\left[\left(1-x/\sqrt{6}\right)^2-\left(1-t/\tilde t\,\right)^{2/3}\right],&\quad 0\leq x \leq \tilde x(t),\\
0,&\quad \tilde x(t)<x\leq 5,
\end{cases}
\end{equation*}
where
$$\tilde{x}(t)=\sqrt{6}(1-(1-t/\tilde{t}\,)^{1/3}),\qquad 0\leq t\leq T<\tilde{t}.$$
This is another non-smooth solution of equation (\ref{heateq}); it exists on the finite time interval $[0,\tilde t\,)$, blowing up when $t$ approaches $\tilde t$.

\begin{table}[t!]
\caption{Errors in solution and conservation laws at the final time $T=10$, when solving {\rm (\ref{heateq})} with conditions in {\rm (\ref{init2})}, using {\rm CS}$(0,\beta)$ methods with $\Delta x=0.25$, $\Delta t=0.333$ and {\rm ML/IM} with $\Delta x=0.25$, $\Delta t=0.04$.}
\small
	\centerline{{\begin{tabular}{ccccc}
\hline
Method &  $\text{Err}_1$ & $\text{Err}_2$ &  Solution error & Computational time\\ 
\hline
$\mbox{CS}(0,-1/4)$	& 3.05e-16   & 1.39e-17 & 0.0162 &  0.013 \\
$\mbox{CS}(0,-1/8)$	&  8.33e-17   & 1.39e-17 & 0.0099 & 0.013\\
$\mbox{CS}(0,0)$	&  2.78e-17   & 1.39e-17 & 0.0036 & 0.002 \\
$\mbox{CS}(0,0.06)$	&  5.55e-17   & 2.08e-17 & 0.0017 & 0.011 \\
ML/IM &  0.0099   & 0.0117 & 0.0165 & 0.076\\ 
\hline
\end{tabular}}}
\label{num3tab1}
\end{table}
\begin{table}[t!]
\caption{Errors in solution and conservation laws at the final time $T=10$, when solving {\rm(\ref{heateq})} with conditions in {\rm(\ref{init2})}, using {\rm CS}$(0,\beta)$ methods with $\Delta x=0.1$, $\Delta t=0.133$ and {\rm ML/IM} with $\Delta x=0.1$, $\Delta t=0.00667$.}
\small
	\centerline{\begin{tabular}{ccccc}
\hline
Method &  $\text{Err}_1$ & $\text{Err}_2$ &  Solution error & Computational time\\ 
\hline
$\mbox{CS}(0,-1/4)$	& 1.78e-16   & 6.66e-17 & 0.0030 &  0.05 \\
$\mbox{CS}(0,-1/8)$	&  2.22e-16   & 4.44e-17 & 0.0018 & 0.04\\
$\mbox{CS}(0,0)$	&  1.33e-16   & 1.33e-16 & 7.01e-04 & 0.01 \\
$\mbox{CS}(0,0.05)$	&  5.33e-16   & 4.44e-17 & 5.23e-04 & 0.04 \\
ML/IM &  0.0036   & 0.0044 & 0.0088 & 0.91\\ 
\hline
\end{tabular}}
\label{num3tab2}
\end{table}

Tables~\ref{num3tab1} and \ref{num3tab2} summarise the numerical solutions of the nonlinear heat equation over the time interval $[0,10]$, with the initial and boundary conditions in (\ref{init2}) for $\tilde{t}=11$. For the coarser grid with $\Delta x=0.25$, the optimal value of $\beta$ is $0.06$; for the finer grid, $\Delta x=0.1$, $\beta=0.05$ is optimal. On both grids, we have chosen $\Delta t>\Delta x$ for the $\mbox{CS}(0,\beta)$ schemes, to show that this does not produce instability. This contrasts markedly with 
ML/IM, which requires $\Delta t<\Delta x^2$ for convergence.
Once again, $\mbox{CS}(0,0)$ shows itself to be a highly-efficient, reasonably-accurate scheme. 
\begin{figure}[htbp]
\centering{\includegraphics[width=14cm,height=17cm]{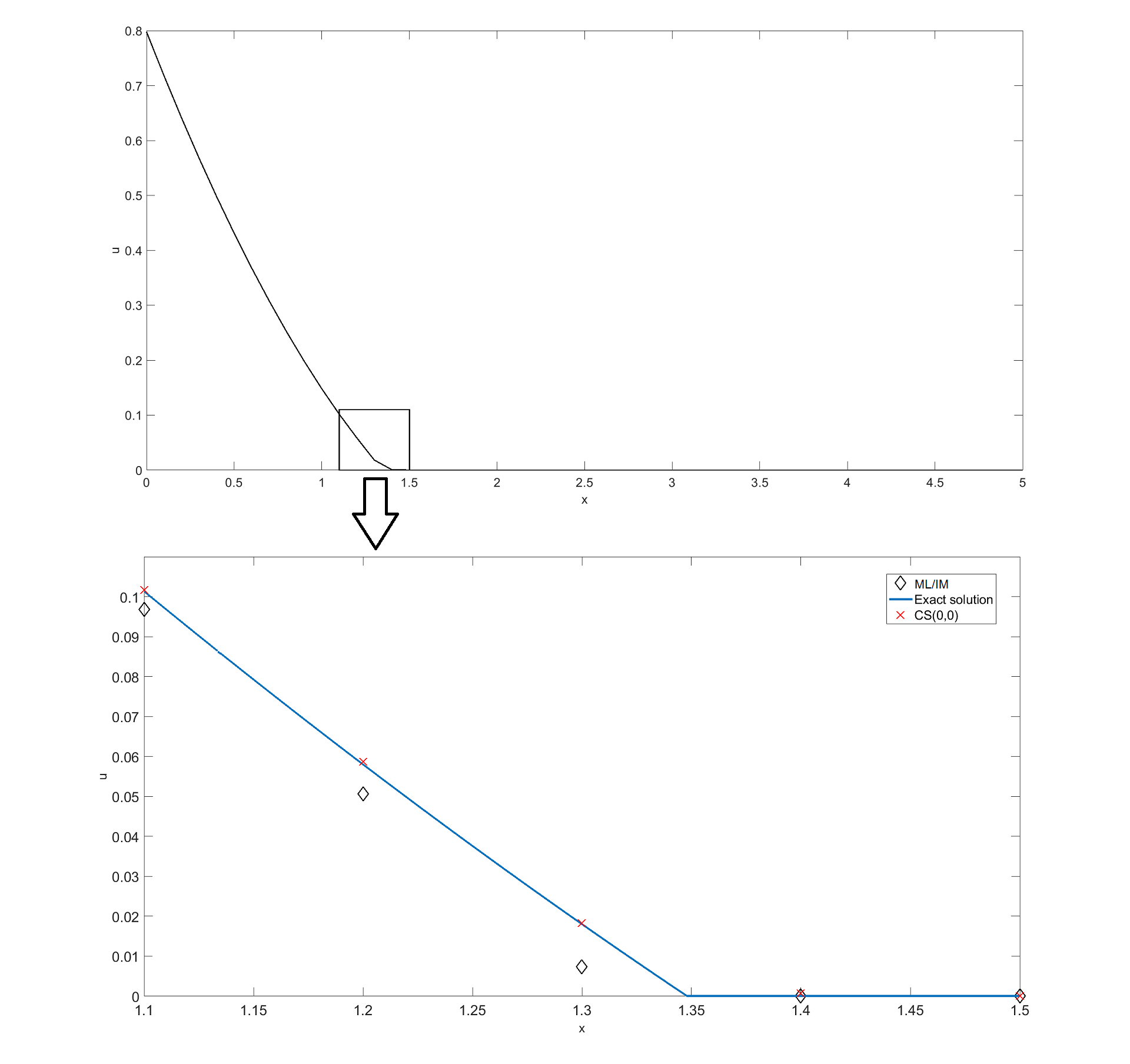}}
\caption{ Results for equation (\ref{heateq}) with conditions in (\ref{init2}). Top: Numerical solution given by method $\mbox{CS}(0,0)$, setting {$\Delta x=0.1$ and $\Delta t=0.133$} over $[0,5]$ at time $T=10$. Bottom: Exact solution (solid curve), numerical solutions close to the interface from $\mbox{CS}(0,0)$ with $\Delta x=0.1$ and $\Delta t=0.133$ (crosses), and ML/IM with $\Delta x=0.1$ and $\Delta t=0.00667$ (diamonds).}
\label{fig2}
\end{figure}  

Figure~\ref{fig2} compares the numerical solutions from $\mbox{CS}(0,0)$ and ML/IM on the finer grid. Again, $\mbox{CS}(0,0)$ is very close to the true solution (as are the other CS schemes in the tables), whereas ML/IM has a small lag close to the interface.

\section{Conclusions and discussion}\label{conclusions}
Motivated by the basic principle of geometric integration that numerical schemes should preserve key structural features of the approximated problem to the extent that is possible, we have presented a strategy for developing finite difference methods that preserve two local conservation laws. This new strategy simplifies the approach introduced in \citet{GrantHydon} and developed in \citet{Grant}. Depending on the stencil, Grant's method can have a very long symbolic computation time (typically several days on a fast PC for a 10-point stencil), which is a strong limitation. However, this difficulty can be overcome by restricting attention to schemes that are second-order, with key terms that are as compact as possible. Such schemes can be determined by hand, or by a short symbolic computation (of no more than a few minutes), even for larger stencils. 

We have developed new parametrized families of conservative numerical schemes for the solution of the KdV equation and a nonlinear heat equation. These schemes seem to be more robust and efficient than other well-known methods that do not preserve multiple conservation laws, perhaps due to topological and analytic advantages. Conservation laws have a topological origin as cohomology classes in the restricted variational bicomplex; schemes that preserve these retain discrete analogues of what may be essential topological features. Furthermore, parameters typically multiply terms that regularize the approximation of $\mca$ in some way. By using benchmark problems and optimising the parameters with respect to the solution error, we have found members of each family that are highly accurate. In practice, the exact solution to a given problem is not usually known. Nevertheless, one can optimise the parameters numerically in order to achieve the best regularization for a given problem. Depending on the problem being approximated, it may be advantageous to choose the parameters in a way that best preserves other geometric structures, such as symplecticity, symmetries, or further conservation laws.

\subsection*{Acknowledgements}
We are grateful to our colleague John Pearson, University of Edinburgh, and to the anonymous referees, for their insightful remarks and constructive suggestions which have helped to improve this paper.

\end{document}